\documentclass[8 pt]{amsproc}
\usepackage{amsmath}
\usepackage{graphicx}
\usepackage{amssymb}
\oddsidemargin=-1cm\evensidemargin=-1cm
\begin{document}
\numberwithin{equation}{section}

\def\1#1{\overline{#1}}
\def\2#1{\widetilde{#1}}
\def\3#1{\widehat{#1}}
\def\4#1{\mathbb{#1}}
\def\5#1{\frak{#1}}
\def\6#1{{\mathcal{#1}}}

\def\C{{\4C}}
\def\R{{\4R}}
\def\N{{\4N}}
\def\Z{{\4Z}}

\title{CR Singularities and Generalizations of Moser's Theorem I}
 
\author{Valentin Burcea}
\begin{abstract}  It is studied the convergence for Formal Holomorphic Mappings defined between Special Classes of C.-R. Singular Real-Analytic Submanifolds in Complex Spaces. In particular,   it is solved the local Equivalence Problem in some particular cases.  \end{abstract}

\address{INDEPENDENT}
\email{vdburcea@gmail.com}

\thanks{ Thanks($20\%$) to CAPES at The Federal University of Minas Gerais, for supporting that visit of $4$ months; Special Thanks($20\%$) to Trinity College Dublin, because I received also a Posgraduate Fellowship in the period $2010-2011$; Special Thanks($20\%$) to Science Foundation Ireland Grant 10/RFP/MT H2878; Special Thanks($40\%$) also to Science Foundation Ireland Grant 06/RFP/MAT 018 on this paper, emphasizing that the reference \cite{V1} was fully supported by  Science Foundation Ireland Grant 06/RFP/MAT 018. }
\maketitle

\def\cn{{\C^n}}
\def\cnn{{\C^{n'}}}
\def\ocn{\2{\C^n}}
\def\ocnn{\2{\C^{n'}}}


\def\dist{{\rm dist}}
\def\const{{\rm const}}
\def\rk{{\rm rank\,}}
\def\id{{\sf id}}
\def\tr{{\bf tr\,}}
\def\aut{{\sf aut}}
\def\Aut{{\sf Aut}}
\def\CR{{\rm CR}}
\def\GL{{\sf GL}}
\def\Re{{\sf Re}\,}
\def\Im{{\sf Im}\,}
\def\span{\text{\rm span}}
\def\Diff{{\sf Diff}}

\def\codim{{\rm codim}}
\def\crd{\dim_{{\rm CR}}}
\def\crc{{\rm codim_{CR}}}

\def\phi{\varphi}
\def\eps{\varepsilon}
\def\d{\partial}
\def\a{\alpha}
\def\b{\beta}
\def\g{\gamma}
\def\G{\Gamma}
\def\D{\Delta}
\def\Om{\Omega}
\def\k{\kappa}
\def\l{\lambda}
\def\L{\Lambda}
\def\z{{\bar z}}
\def\w{{\bar w}}
\def\Z{{\1Z}}
\def\t{\tau}
\def\th{\theta}

\emergencystretch15pt \frenchspacing

\newtheorem{Thm}{Theorem}[section]
\newtheorem{Cor}[Thm]{Corollary}
\newtheorem{Pro}[Thm]{Proposition}
\newtheorem{Lem}[Thm]{Lemma}

\theoremstyle{definition}\newtheorem{Def}[Thm]{Definition}

\theoremstyle{remark}
\newtheorem{Rem}[Thm]{Remark}
\newtheorem{Exa}[Thm]{Example}
\newtheorem{Exs}[Thm]{Examples}

\def\bl{\begin{Lem}}
\def\el{\end{Lem}}
\def\bp{\begin{Pro}}
\def\ep{\end{Pro}}
\def\bt{\begin{Thm}}
\def\et{\end{Thm}}
\def\bc{\begin{Cor}}
\def\ec{\end{Cor}}
\def\bd{\begin{Def}}
\def\ed{\end{Def}}
\def\br{\begin{Rem}}
\def\er{\end{Rem}}
\def\be{\begin{Exa}}
\def\ee{\end{Exa}}
\def\bpf{\begin{proof}}
\def\epf{\end{proof}}
\def\ben{\begin{enumerate}}
\def\een{\end{enumerate}}
\def\beq{\begin{equation}}
\def\eeq{\end{equation}}
\section{Introduction and Main Results}

 Two  Real-Analytic Submanifolds  in Complex Spaces, may be only formally (holomorphically) equivalent   like it has been shown by Moser-Webster\cite{MW} and  Gong\cite{G2} in the C.-R. Singular Situation\cite{Bi}, and respectively by  Kossovskiy-Shafikov\cite{RI} in the  C.-R. Situation\cite{BREbook}. Such aspect motivates the local Equivalence Problem \cite{ko2},\cite{ko3} in Complex Analysis, which asks when two Real-Analytic formally equivalent Submanifolds   are actually holomorphically equivalent in Complex Spaces.  

 Moser\cite{Mos} considered the   Equivalence Problem for the Class of  Real-Analytic Surfaces in $\mathbb{C}^{2}$ defined near $p=0$ by
\begin{equation}w=z\overline{z}+\mbox{O}\left(\left|z\right|^{3}\right),\label{m11}\end{equation} 
where $(z,w)$ are the coordinates from $\mathbb{C}^{2}$, showing that if (\ref{m11}) is formally equivalent to the Model 
\begin{equation}w=z\overline{z}\label{m111},\end{equation} 
then (\ref{m11}) is holomorphically equivalent to (\ref{m111}). This result is widely known as the   Theorem of Moser\cite{Mos}.   

 In this paper, we work with large Classes of Real-Smooth Submanifolds  of Codimension $2d$ in $\mathbb{C}^{N+d}$, defined near $p=0\in \mathbb{C}^{N+d}$,
 in order to extend the Theorem of Moser\cite{Mos} from $\mathbb{C}^{2}$ to  $\mathbb{C}^{N+d}$, in  the coordinates     \begin{equation}(z;w):=\left(z_{1},z_{2},\dots,z_{N};w_{1},w_{2},\dots,w_{d}\right)\in \mathbb{C}^{N+d},\quad\mbox{where $N,d\geq 1$.}\label{coord}\end{equation}

In particular, we consider $M\subset\mathbb{C}^{N+d}$   a Real-Smooth Submanifold defined  near   $p=0\in \mathbb{C}^{N+d}$ by
\begin{equation}w=Q(z,\overline{z})+\mbox{O}\left(\left|z\right|^{3}\right), \label{M1}\end{equation}
where  the vector polynomial   of degree $2$ in $(z,\overline{z})$, arising from the  formal expansion of its local defining equation, denoted by \begin{equation}Q(z,\overline{z})=\left(Q_{1}(z,\overline{z}),Q_{2}(z,\overline{z}),\dots,Q_{d}(z,\overline{z})\right).\label{Qu}\end{equation}  

The Real-Submanifold $M$, defined by (\ref{M1}), is called   Non-Degenerate if   
 \begin{equation}\frac{\partial Q}{\partial z_{1}}(z,\overline{z})= 0,\frac{\partial Q}{\partial z_{2}}(z,\overline{z})= 0,\dots,\frac{\partial Q}{\partial z_{N}}(z,\overline{z})= 0 \mapsto  z=0\in\mathbb{C}^{N}.\label{Deg}
  \end{equation}  

The Real-Submanifold $M$, defined by (\ref{M1}), is called   Semi-Real if  
\begin{equation}Q_{1}(z,\overline{z})=\overline{Q_{1}(z,\overline{z})},Q_{2}(z,\overline{z})=\overline{Q_{2}(z,\overline{z})},\dots,Q_{d}(z,\overline{z})=\overline{Q_{d}(z,\overline{z})}.\label{real}\end{equation}

The Real-Submanifold $M$, defined by (\ref{M1}), is called   $d$-ranked if 
\begin{equation}\mbox{Rank}\begin{pmatrix}\frac{\partial Q_{1}}{\partial z_{1}}(z,\overline{z})& \frac{\partial  Q_{1}}{\partial z_{2}}(z,\overline{z}) & \dots &  \frac{\partial Q_{1}}{\partial z_{N}}(z,\overline{z})
\\  \frac{\partial Q_{2}}{\partial z_{1}}(z,\overline{z})& \frac{\partial Q_{2}}{\partial z_{2}}(z,\overline{z}) & \dots &  \frac{\partial Q_{2}}{\partial z_{N}}(z,\overline{z}) \\  \vdots & \vdots & \ddots &  \vdots \\  \frac{\partial Q_{d}}{\partial z_{1}}(z,\overline{z})& \frac{\partial Q_{d}}{\partial z_{2}}(z,\overline{z}) & \dots &  \frac{\partial Q_{d}}{\partial z_{N}}(z,\overline{z})
\end{pmatrix}=d.\label{Degg}
\end{equation}

The point  $p=0\in\mathbb{C}^{N+d}$ is a C.-R. Singularity\cite{V1},\cite{V11}, for the Real-Submanifold $M$ defined by (\ref{M1}), if it is
a jumping point for the mapping  $M\ni p \mapsto \dim_{\mathbb{R}}\left(T_{p}^{c}M  \right)$  defined near $p=0\in M\subset\mathbb{C}^{N+d}$. The conditions (\ref{Deg})   and  (\ref{Degg}) define independent Classes of Real-Smooth Submanifolds, if (\ref{real}) holds.  In particular, we obtain

 \bt\label{t1}  Let $M,M'\subset\mathbb{C}^{N+d}$ be Semi-Real Non-Degenerate Real-Analytic Submanifolds, provided that $M$ is $d$-ranked and $N>d>1$.  Then,  any formal mapping, from $M$  into $M'$, is a convergent equivalence.  \et 
  
 This result relies on a  careful analysis of the considered Formal (Holomorphic) Mapping  in the local  defining equation   of $M'$  when    (\ref{M1}), (\ref{real}) and (\ref{Degg})  hold. In particular,  we  use a family of   Real-Analytic Submanifolds in  $\mathbb{C}^{N}$, regardless that such  Submanifolds  may not be minimal. They   exist  near   $p=0\in \mathbb{C}^{N+d}$ because (\ref{Degg})  holds and $N>d$.  In particular, we obtain
  \bc\label{ct1}  Let $M,M'\subset\mathbb{C}^{N+d}$ be Semi-Real NonDegenerate Real-Smooth Submanifolds, provided $N>d>1$ and that $M$ is $d$-ranked.  Then,  any non-constant formal mapping, from $M$  into $M'$, is an   equivalence.  \ec 
 
Standardly  considered   Baouendi-Mir-Rothschild\cite{BMR1}, Mir\cite{Mir1},\cite{Mir2}, Suny\'{e}\cite{Su}, Merker\cite{Mer}, Meylan-Mir-Zaitsev\cite{MNZ},   the minimality  is   an optimal geometrical condition  in order to obtain the convergence  of Formal Holomorphic Mappings  between Real-Analytic Submanifolds in Complex Spaces. Its absence may generate    interesting phenomenons.   Kossovskiy-Shafikov\cite{RI} showed recently that it may not  exist   Holomorphic Equivalences between two non-minimal Real-Analytic  Formally  Equivalent  Submanifolds in   Complex Spaces.   

However, the minimality is less relevant  through this paper, regardless  that a Real-Submanifold may not be minimal near its C.-R. Singularities\cite{Dol1},\cite{Dol3},\cite{DTZ1},\cite{DTZ2}.  It suffices to prove the convergence of the Formal Mappings at points near   $p=0\in \mathbb{C}^{N+d}$, because their Holomorphicity is provided by The Phenomenon of Hartogs. The local convergences are provided by 
 Artin's Approximation Theorem\cite{A} combined with  the  convergence of the formal mapping restricted to their first Segre mappings\cite{BREbook},\cite{Mir1}. In particular, the assumption $d>1$ suffices  in order 
to construct convenient analytic systems. The case $d=1$ is different, because  it  does not provide generally enough analytic equations in order to construct convenient analytic systems. It is required to consider other assumptions like:

\bt\label{t2}  Let $M,M'\subset\mathbb{C}^{N+1}$ be  Real-Analytic Semi-Real Non-Degenerate Submanifolds such that    $M'$ can not   be formally   transformed into the   Model
  \begin{equation}w=Q(z,\overline{z}) .\label{Model}\end{equation}
  
 Then provided   $N>1$, any formal mapping from $M$ is into $M'$, is a convergent equivalence.
 \et  
 
The  proof of Theorem \ref{t2} follows  the proof of Theorem \ref{t1}. Because $M'$ can not   be formally   transformed into  Model (\ref{Model}), it follows the existence of  terms of degree at least $3$  in the local defining equation of $M'$. It  may be seen  like a  non-degeneracy condition in order to derive    analytic systems from its local defining equation. Their absence gives an insufficient number of  analytic equations, because the polynomial $Q(z,\overline{z})$ is real-valued. In particular, we  obtain

 \bc\label{ct2}  Let $M,M'\subset\mathbb{C}^{N+1}$ be a Real-Smooth Semi-Real Non-Degenerate Submanifolds such that    $M$ can not   be formally   transformed into     Model (\ref{Model}).  Then provided   $N>1$, any non-constant formal mapping from $M$ is into $M'$, is an equivalence.
 \ec  
 
 It is interesting to observe the following:
 \bc Let $M,M'\subset\mathbb{C}^{N+1}$ be formally equivalent Real-Analytic Semi-Real NonDegenerate Submanifolds such that    $M$ can not   be formally   transformed into the   Model
  \begin{equation}w=z_{1}\overline{z}_{1}+z_{2}\overline{z}_{2}+\dots+z_{N}\overline{z}_{N}+ \lambda_{1}\left(z_{1}^{2}+\overline{z}_{1}^{2}\right)+\lambda_{2}\left(z_{2}^{2}+\overline{z}_{2}^{2}\right)+\dots+\lambda_{N}\left(z_{N}^{2}+\overline{z}_{N}^{2}\right).\label{Modell}\end{equation}
  
 Then, $M$ and $M'$ are biholomorphically equivalent,  provided that  $N>1$ and $\lambda_{1},\lambda_{2},\dots,\lambda_{N}\geq 0$.
 \ec 
 
Such  C.-R. Singularities are called Special C.-R. Singularities\cite{Dol1},\cite{Dol3}. These real numbers $\lambda_{1},\lambda_{2},\dots,\lambda_{N}$, are called Generalized Bishop invariants and, are important in order to study  local Equivalence Problems.  In $\mathbb{C}^{2}$, Gong\cite{G1} proved another analogue of Moser's Theorem\cite{Mos} when the Bishop invariant\cite{Bi} is not vanishing.  For $\lambda_{1}=\lambda_{2}=\dots=\lambda_{N}=0$, Huang-Yin\cite{HY1}  generalized the  Theorem of Moser\cite{Mos}. They\cite{HY1}  proved that  (\ref{M1}) is formally equivalent to the Model  (\ref{Modell})
if and only if their holomorphic equivalence occurs.  

Next, we consider $M\subset\mathbb{C}^{N+1}$   a Real-Smooth Submanifold defined  near   $p=0\in \mathbb{C}^{N+d}$ by
\begin{equation}w=Q(z,\overline{z})+\mbox{O}\left(\left|z\right|^{3}\right), \label{M2}\end{equation}
where  the  polynomial   of degree $2$ in $(z,\overline{z})$,   arising from the  formal expansion of its local defining equation, is denoted by $Q(z,\overline{z})$. 

The Real-Submanifold $M$, defined by (\ref{M2}), is called Semi-Complex if  
\begin{equation}Q(z,\overline{z})\not\equiv\overline{Q(z,\overline{z})}.\label{real1}\end{equation}

The Real-Submanifold $M$, defined by (\ref{M2}), is called  ranked if 
\begin{equation}\mbox{Rank}\begin{pmatrix}\frac{\partial \left(\Re Q\right)}{\partial z_{1}}(z,\overline{z})& \frac{\partial \left(\Re Q\right)}{\partial z_{2}}(z,\overline{z}) & \dots &  \frac{\partial \left(\Re Q\right)}{\partial z_{N}}(z,\overline{z})
\\  \frac{\partial \left(\Im Q\right)}{\partial z_{1}}(z,\overline{z})& \frac{\partial \left(\Im Q\right)}{\partial z_{2}}(z,\overline{z}) & \dots &  \frac{\partial \left(\Im Q\right)}{\partial z_{N}}(z,\overline{z})
\end{pmatrix}=2.\label{Degg1}  
\end{equation}

The conditions (\ref{Deg})   and  (\ref{Degg1}) define independent Classes of Real-Smooth Submanifolds, if (\ref{real1}) holds.  The  polynomial of degree $2$ denoted by $Q(z,\overline{z})$ is generally complex-valued (see \cite{V11}), excepting when   $N=1$.  In particular, we obtain

\bt\label{t3}  Let $M,M'\subset\mathbb{C}^{N+1}$ be Semi-Complex Non-Degenerate Real-Analytic Submanifolds provided that $M$ is  ranked and $N>2$.  Then,  any formal mapping, from $M$  into $M'$, is a convergent equivalence.  \et 

 This result relies on a  careful analysis of the considered Formal (Holomorphic) Mapping  in the local  defining equation   of $M'$  when    (\ref{M2}), (\ref{real1}) and (\ref{Degg1})  hold. We focus on a family of possibly non-minimal  Real-Analytic Submanifolds in  $\mathbb{C}^{N}$. They   exist  near   $p=0\in \mathbb{C}^{N+d}$ because (\ref{Degg1})  holds and $N>2$.  It suffices to prove the convergence of the Formal Mappings at points near   $p=0\in \mathbb{C}^{N+d}$, because their Holomorphicity is provided by The Phenomenon of Hartogs. The local convergences are provided by 
 Artin's Approximation Theorem\cite{A} combined with  the  convergence of the formal mapping restricted to their first Segre mappings\cite{BREbook},\cite{Mir1}.  The non-reality (\ref{real1}) of  the polynomial $Q(z,\overline{z})$ is   optimal   in order to derive by (\ref{Degg1}) convenient  analytic systems 
 from (\ref{M2}), because its terms of degree $3$ are not relevant. In particular, we obtain

\bc\label{ct3}  Let $M,M'\subset\mathbb{C}^{N+1}$ be Semi-Complex Non-Degenerate Real-Smooth Submanifolds provided that $M$ is  ranked and $N>2$.  Then,  any non-constant formal mapping, from $M$  into $M'$,
  is an equivalence.  \ec 

It is an interesting property of rigidity for formal mappings, because  the existence of a formal mapping, between two Real-Smooth Submanifolds in Complex Spaces, does not imply  the existence of a formal equivalence between two Real-Smooth Submanifolds in Complex Spaces. The Smoothness of Real Submanifolds is essential in order to attain this property in the previous assumptions.  Generally, it is possible to deal  with formal mappings that are not equivalences. Their convergence is obtained  by reiterating the strategy of the proofs of Theorems \ref{t1}, \ref{t2} and \ref{t3} in  previous circumstances. In particular, we obtain 
 \bt\label{t4}  Let $M,M'\subset\mathbb{C}^{N+1}$ be  Real-Analytic  Non-Degenerate Submanifolds of codimension $2$ such that    $M'$ is Semi-Real and can not   be formally   transformed into
   (\ref{Model}), and $M$ is Semi-Complex.  Then provided   $N>1$, any formal mapping from $M$   into $M'$, is convergent.
 \et  
 
The  proof of Theorem \ref{t4} uses the lines of  the proof of Theorem \ref{t2}, because
the existence of   terms of degree at least $3$, in its local defining equation,   may be seen  like an optimal  non-degeneracy condition in order to construct convenient  analytic systems. Their absence gives an insufficient number of  analytic equations, because the polynomial $Q(z,\overline{z})$ is real-valued. 

On the other hand, if the polynomial $Q(z,\overline{z})$ is not real-valued, we obtain 
\bt\label{t5}  Let $M,M'\subset\mathbb{C}^{N+1}$ be a Real-Analytic  Non-Degenerate Submanifolds of codimension $2$, provided  that    $M'$ is Semi-Complex and Ranked  and $M$ is SemiReal.  Then provided   $N>2$, any formal mapping from $M$   into $M'$, is convergent.
 \et  
 
The  proof of Theorem \ref{t5} uses the lines of  the proof of Theorem \ref{t2}, because
the non-reality of the polynomial $Q(z,\overline{z})$  may be seen  like an optimal  non-degeneracy condition in order to construct convenient  analytic systems, but replacing the target Real-Submanifold  of codimension $2$ with a Real Submanifold assumed  Semi-Real and $d$-ranked in $\mathbb{C}^{n+d}$ with $d\geq 2$, it exists an appropriate background in order to construct convenient analytic systems following the lines of Theorem \ref{t1}. In particular, we obtain

 \bt\label{t6}  Let $M,M'\subset\mathbb{C}^{N+1}$ be a Real-Analytic  Non-Degenerate Submanifolds, provided  that    $M'$ is Semi-Complex, Ranked and of codimension at $2$ and $M$ is Semi-Real and $d$-ranked with $d\geq 2$.  Then provided   $N>2$, any formal mapping from $M'$   into $M$, is convergent.
 \et  
 
The constructions of analytic systems are pursued 
according to  the guidance, received by the author\cite{V1} during the his initiation\cite{V1},  from  Zaitsev\cite{D11},\cite{D1} with  appreciations also towards to his  doctoral learnings\cite{V11}. The analytic equations, which form  analytic systems, are obtained by identifying the coefficients of the powers of   parameters from local defining equations. This approach is motivated by the standard procedure  considered usually in order to construct  normal forms\cite{ko2},\cite{ko3},\cite{BEL5},\cite{D11},\cite{D1}, providing  crucial arguments through   this paper. 

We introduce the Class of Real-Analytic Hypersurfaces 
\begin{equation}\mathbb{C}^{N'+1} \supset M':\hspace{0.1 cm}
\Im w'= \tilde{\varphi}\left(z',\overline{z}',\Re w'\right),\label{M3}
\end{equation}
using to the following formal expansion
\begin{equation}\tilde{\varphi}\left(z',\overline{z}',\Re w'\right)=\displaystyle\sum_{m,k\in\mathbb{N}\atop{m+k \geq 2}}\left(\Re w'\right)^{m} \cdot P_{m\hspace{0.02 cm} k}\left(z',\overline{z}'\right),\quad\mbox{ defined near}\hspace{0.1 cm}p=0\in\mathbb{C}^{N'+1},\label{pee1}
\end{equation}
where we have used suitable homogeneous polynomials in $(z,\overline{z})$ in the coordinates
$\left(w';z'\right):=\left(w';z'_{1},z'_{2},\dots,z'_{N'}\right)\in\mathbb{C}^{N'+1}$, for $N'\in\mathbb{N}^{\star}$.

The Real-Analytic Hypersurface $M'$ is called NonDegenerate if the polynomial $P_{m_{0}\hspace{0.02 cm} k_{0}}$  satisfies (\ref{Deg}),  where 
\begin{equation} \left(m_{0},k_{0}\right)=\mbox{Min}\left\{\left(m,k\right)\in\mathbb{N}\times \mathbb{N};\hspace{0.1 cm}P_{m  \hspace{0.02 cm} k}\not\equiv 0\right\}.\label{min}
\end{equation}

Going  forward, the  hypothesises (\ref{M1}) and (\ref{M2})  may be used together with their non-degeneracy conditions (\ref{Qu}),(\ref{Deg}),(\ref{Degg})  and (\ref{Degg1}) with (\ref{real}) and (\ref{real1}) in order to prove other results of convergence.  In particular, we obtain

\bt \label{t7}Let $M\subset\mathbb{C}^{N+1}$ be a Semi-Real Non-Degenerate Real-Analytic Submanifold of Codimension $2$. Then, any formal mapping, from $M$ into $M'$, is convergent, provided  that  $M'\subset \mathbb{C}^{N'+1}$ is a Real-Analytic  NonDegenerate Hypersurface and    $N>1$.    
\et

 The NonDegeneracy of $M'$ seems to be   optimal   in order to  pursue (formal) computations of normal form type, regardless that we deal with formal mappings which are not equivalences, but we can make certain non-degeneracy assumption on the formal mappings. Like in the proof of Theorem \ref{t1}, we focus on a family of eventually non-minimal  Real-Analytic Submanifolds in  $\mathbb{C}^{N}$.  In particular, we prove  the partial convergence of the formal mapping restricted to their first Segre mappings\cite{BREbook},\cite{Mir1}, regardless that $M'$ may not be minimal  because $m_{0}$ may not vanish. Then     Artin's Approximation Theorem applies in order to prove the convergence of the formal mapping at points near $p=0\in\mathbb{C}^{N+d}$.  Its Holomorphicity is provided by the Phenomenon of Hartogs.  In particular, we obtain 
 \bt \label{t8}Let $M\subset\mathbb{C}^{N+1}$ be a Ranked Semi-Complex Non-Degenerate Real-Analytic    Submanifold of Codimension $2$. Then, any formal mapping, from $M$ into $M'$, is convergent, provided that  $M'\subset \mathbb{C}^{N'+1}$ is a Real-Analytic   Non-Degenerate Hypersurface  and   $N>2$.  
\et

In particular, it follows that the result of the initial version of this work holds. Its mistakes have been an inspiration source for me   in order to develop this paper by multiplying for $10$ times the initial result. The Financial Thanks have been re-organized like the relevance of each Funding, emphasizing that  I did not receive the first monthly scholarship ($4100$ BRL) in Minas Gerais. Contrary, I recall that I received    a lost sum  of money from The School of Mathematics of Trinity College Dublin, and again also   from Science Foundation of Ireland, before moving to Brazil. I felt joy in each of my  visits in various cities from Brazil, but the mistakes effectuated continue to provide inspiration. In particular, I will never give academic presentations in Brazilian Universities, because there exist unrepairable prejudices. 

Finally, we introduce the Special Class of Real-Analytic Surfaces 
\begin{equation}\mathbb{C}^{2} \supset M': 
 \hspace{0.1 cm}  w'=P\left(z',\overline{z}'\right)+\mbox{O}\left(k_{0}+1\right) ,\quad \mbox{where   $P\not\equiv 0$ and $\mbox{Deg}\left(P\right)=k_{0}$,} \label{M5}
\end{equation}
 defined near   $p=0\in\mathbb{C}^{2}$  in the coordinates (\ref{coord}) for $d=1$ and $N=1$.
It defines a Special Class of Degenerate C.-R. Singularities.  The degree of the homogeneous
 polynomial $P$, denoted by $k_{0}$, is assumed to be at least $2$. The polynomial $P$ may be assumed real-valued\cite{Bi} only when $k_{0}=2$. The polynomial  $P$ may not be real-valued\cite{Bha1},\cite{Bha2},\cite{Bha3} when $k_{0}>2$. It is an interesting case that deserves to be studied. The interactions of terms become complicated in the local defining equations especially when $k_{0}>2$. In particular, we obtain

\bt \label{t9}Let $M\subset\mathbb{C}^{N+1}$ be a Semi-Real Non-Degenerate Real-Analytic Submanifold of Codimension $2$. Then, any formal mapping, from $M$ into $M'$, is convergent, provided  that  $M'\subset \mathbb{C}^{2}$ is a Real-Analytic  Surface defined as in (\ref{M5}).    
\et

 The proof of Theorem \ref{t9} is very related to the proof of Theorem \ref{t2} when $k_{0}=2$.  The proof  shows non-trivialities when $k_{0}>2$, because the possible non-reality of the polynomial $P$ creates computational difficulties.  In particular, we obtain
\bt \label{t10}Let $M\subset\mathbb{C}^{N+1}$ be a Semi-Complex Ranked Non-Degenerate Real-Analytic Submanifold of Codimension $2$. Then, any formal mapping, from $M$ into $M'$, is convergent, provided  that  $M'\subset \mathbb{C}^{2}$ is a Real-Analytic  Surface defined as in (\ref{M5}).    
\et

The  C.-R. Singularities \cite{Bi},\cite{V11},\cite{DTZ1},\cite{GL},\cite{HK},\cite{HY2221},\cite{Sla2}   form   an important subject for Complex Analysis. Surfaces in $\mathbb{C}^{2}$ with degenerate CR singularities and  local polynomial convexity problems near degenerate CR singularities have
been studied by Bharali\cite{Bha1},\cite{Bha2},\cite{Bha3}. Gong-Lebl\cite{GL}   studied    C.-R. Singular Submanifolds of codimension $2$ which are Levi-flat at the C.-R. points. Slapar \cite{Sla2} showed that a compact Real Surface embedded in a Complex Surface has a regular Stein neighborhood basis assuming the existence of C.-R. Singularities on the Real Surface. Dolbeault\cite{Dol1}, \cite{Dol3}, Dolbeault-Tomassini-Zaitsev\cite{DTZ2},\cite{DTZ1}   used the existence of the C.-R. Singularities  in order to study   problems as the  existence and the uniqueness of the Levi-flat hypersurfaces with prescribed compact boundary\cite{DTZ2},\cite{DTZ1}. The author\cite{V11} constructed a family of analytic discs attached to a Class of  C.-R. Singular Real Submanifolds in codimension $2$  trying to understand the local hull of holomorphy using methods from Huang-Krantz\cite{HK}.  Furthermore, Huang-Yin\cite{HY2221},\cite{HY2222}  impressively exploited the C.-R. structure near the C.-R. singularity\cite{Bi},\cite{V11} in order to study the local hull of holomorphy\cite{HY2221},\cite{HY2222}.

Few words about the organization of this paper. In Section $2$, there are discussed the non-degeneracy conditions (see (\ref{Deg}), (\ref{Degg}), (\ref{Degg1})) and the notations considered throughout this paper. In  Sections $3$ and $4$, there are discussed the used families of Real-Analytic Submanifolds, starting with Section $5$, in order to provide the arguments of convergence in the final Sections of this paper. 

\subsection{Acknowledgements} Special Thanks to the Editorial Board of   \emph{Indian Journal of Pure and Applied Mathematics}   for stopping the publication of a flawed version of this paper.   It originates from my period from The Federal University of Minas Gerais. Meanwhile, I developed the initial work by  learning from my own mistakes, also by using learnings\cite{V1},\cite{V11} from my doctoral period. Special Thanks to my supervisor Prof. Dmitri Zaitsev  for bringing \cite{GL} into my attention.  Thanks also to Karen O'Doherty, reiterating that  the main part of my doctoral thesis \cite{V1} was fully supported by Science Foundation Ireland Grant 06/RFP/MAT 018.

\section{Environment} 

\subsection{SemiReal Real-Analytic Submanifolds} Let   
$M$ and $M'$ be  Real-Analytic Submanifolds    near $p=0\in\mathbb{C}^{N+d}$ in the coordinates (\ref{coord}), formally   equivalent and  defined by
\begin{equation}\begin{split}&\hspace{0.1 cm}\mathbb{C}^{N+d}\supset M:\hspace{0.1 cm} w= Q(z,\overline{z})+\varphi(z,\overline{z}),\\& \mathbb{C}^{N+d}\supset M':\hspace{0.1 cm} w= Q'(z,\overline{z})+\varphi'(z,\overline{z}), \end{split} \label{W}
\end{equation}
where we have used   vector-polynomials of degree $2$ in $(z,\overline{z})$, denoted by 
\begin{equation} \begin{split}&\hspace{0.1 cm} Q(z,\overline{z})= \left(
Q_{1}(z,\overline{z}) , Q_{2}(z,\overline{z}) , \dots , Q_{d}(z,\overline{z}) 
\right), \\&  Q'(z,\overline{z})= \left(
Q'_{1}(z,\overline{z}) , Q'_{2}(z,\overline{z}) , \dots , Q'_{d}(z,\overline{z}) 
\right),\end{split}\label{W1}
\end{equation}
and   satisfying (\ref{Deg})  and (\ref{real}), for
\begin{equation}\begin{split}&\hspace{0.1 cm}\varphi(z,\overline{z})=\left(
\varphi_{1}(z,\overline{z}) , \varphi_{2}(z,\overline{z}) , \dots , \varphi_{N}(z,\overline{z})
\right)=\left(
\mbox{O}\left(3\right) , \mbox{O}\left(3\right) , \dots, \mbox{O}\left(3\right) 
\right)=\mbox{O}\left(3\right) ,\\&  \varphi'(z,\overline{z})=\left(
\varphi'_{1}(z,\overline{z}) , \varphi'_{2}(z,\overline{z}) , \dots , \varphi'_{N}(z,\overline{z})
\right)=\left(
\mbox{O}\left(3\right) , \mbox{O}\left(3\right) , \dots, \mbox{O}\left(3\right) 
\right)=\mbox{O}\left(3\right) .\end{split} \label{W2} 
\end{equation}
 
 In order to understand better   (\ref{Deg}) and (\ref{real}), we write  
\begin{equation}\begin{split}& Q_{k}(z,\overline{z})=\displaystyle\sum_{m,n=1}^{N}\left(\overline{a}_{mn}^{\left(k\right)}\overline{z}_{m}\overline{z}_{n}+b_{mn}^{\left(k\right)}z_{m}\overline{z}_{n}+a_{mn}^{\left(k\right)}z_{m} z_{n}\right),\quad\mbox{for all $k=1,\dots,d$,}\\&\hspace{0.16 cm}
 Q(z,\overline{z})=\displaystyle\sum_{m,n=1}^{N}\left(\overline{a}_{mn}\overline{z}_{m}\overline{z}_{n}+b_{mn}z_{m}\overline{z}_{n}+a_{mn}z_{m} z_{n}\right),\quad\hspace{0.1 cm}\mbox{for $d=1$.}\end{split}
\label{LL1}
  \end{equation}
 
 In order to understand better   (\ref{Deg}) and (\ref{real}), we compute 
\begin{equation}\begin{split}&  \frac{\partial   Q_{k}}{\partial z_{l}}(z,\overline{z})= \displaystyle\sum_{ n=1}^{N}\left(b_{nl}^{(k)} \overline{z}_{l}+\left( a_{nl}^{(k)}+ a_{nl}^{(k)}\right)  z_{l}\right),\quad\mbox{for all   $l=1,\dots,N$ and $k=1,\dots,d$,}\\&\hspace{0.16 cm}  \frac{\partial   Q}{\partial z_{l}}(z,\overline{z})= \displaystyle\sum_{ n=1}^{N}\left(b_{nl}  \overline{z}_{l}+\left( a_{nl} + a_{ln} \right)  z_{l}\right),\quad\hspace{0.28 cm}\quad\mbox{for all   $l=1,\dots,N$ and $d=1$.}\end{split}   \label{LL2}
  \end{equation}
  
In particular, (\ref{Deg}) holds  in the following cases particular and independent cases 
 \begin{equation}\det\left(\left(b_{mn}\right)_{1\leq m,n\leq N}\right)\neq 0, \quad
\det\left(\left(a_{mn}+a_{nm}\right)_{1\leq m,n\leq N}\right)\neq 0,  \label{lolll}
\end{equation}  
describing invariant non-degeneracy conditions under any linear   or  
 linear-hermitian change of coordinates.

Generally, other cases may occur like
\begin{equation}\left(b_{mn}\right)_{1\leq m,n\leq N}=\begin{pmatrix} 0& 0& \dots & 0 \\ 0& 1& \dots & 0 \\ \vdots & \vdots & \ddots & \vdots  \\ 0 & 0& \dots & 1 
\end{pmatrix},\quad \left(a_{mn}\right)_{1\leq m,n\leq N}=\begin{pmatrix} 1& 0& \dots & 0 \\ 0& 0& \dots & 0 \\ \vdots & \vdots & \ddots & \vdots  \\ 0 & 0& \dots & 0 
\end{pmatrix},\label{lolllt}
\end{equation}
or other derived cases from (\ref{lolllt}) by linear changes of coordinates, which clearly do not satisfy (\ref{lolll}). 

In particular for $d=1$, (\ref{Degg}) reduces to
\begin{equation}\mbox{Rank}\begin{pmatrix}\frac{\partial Q }{\partial z_{1}}(z,\overline{z})& \frac{\partial  Q }{\partial z_{2}}(z,\overline{z}) & \dots &  \frac{\partial Q }{\partial z_{N}}(z,\overline{z})  \end{pmatrix}=1.\label{Dega}
\end{equation}

We observe that   (\ref{Deg}) implies (\ref{Degg}) when $d=1$, but  when
$Q_{0}(z,\overline{z})=\left(Q_{1}(z,\overline{z}),0,\dots,0\right), 
$ where $Q_{1}(z,\overline{z})$ satisfies (\ref{Dega}), it follows that (\ref{Deg}) does not imply (\ref{Degg}) for any $d\geq 2$. When $Q_{0}(z,\overline{z})=\left(z_{1}\overline{z}_{1} ,z_{2}\overline{z}_{2},\dots,z_{d}\overline{z}_{d},0,\dots,0\right)$ it follows that (\ref{Degg}) does not imply (\ref{Deg}) for any $d\geq 1$. In particular for  $d=1$, we work by (\ref{Deg}) and (\ref{Degg1}) with (\ref{W}), (\ref{W1}) and (\ref{W2})  when (\ref{real1})  holds:
\subsection{SemiComplex Real-Analytic Submanifolds}  Because   by (\ref{real1})   the polynomial $Q(z,\overline{z})$  is not be real-valued, we write 
\begin{equation} Q(z,\overline{z})=\displaystyle\sum_{m,n=1}^{N}\left( c_{mn}\overline{z}_{m}\overline{z}_{n}+b_{mn}z_{m}\overline{z}_{n}+a_{mn}z_{m} z_{n}\right).
\label{qu}
  \end{equation}

When
$Q_{0}(z,\overline{z})=z_{1}\overline{z}_{1}+z_{2}\overline{z}_{2}+z_{1}\overline{z}_{2}+0\cdot z_{3}\overline{z}_{3}+\dots+0\cdot z_{N}\overline{z}_{N}$, it follows  that (\ref{Degg1}) does not imply (\ref{Deg})  when (\ref{real1}) holds, and when $Q_{0}(z,\overline{z})=\left(1+\sqrt{-1}\right)z_{1}\overline{z}_{1}+z_{2}\overline{z}_{2} +\dots+z_{N}\overline{z}_{N}+z_{1}\overline{z}_{2}^{2}$
it follows that  (\ref{Deg}) does not imply (\ref{Degg1}) when (\ref{real1}) holds, 
and    satisfies  also the hypothesis of Theorem \ref{t2}. In this regard, other examples may be derived  by perturbating the quadratic model  with higher order terms, in order to    
 make computations using:
\subsection{Formal (Holomorphic) Mappings} We work with formal (holomorphic) equivalences denoted by \begin{equation}\left(g(z,w);f(z,w)\right)=\left(g_{1}(z,w),g_{2}(z,w),\dots, g_{d}(z,w);f_{1}(z,w),f_{2}(z,w),\dots, f_{N'}(z,w)\right),\quad\mbox{for $d,N '\in\mathbb{N}^{\star}$. }\label{equiv}\end{equation}
   
In particular,  we work with formal expansions 
\begin{equation}\begin{split}& g(z,w)=\displaystyle\sum_{J\in\mathbb{N}^{d}}g_{J}(z)w^{J}=\left(\displaystyle\sum_{J\in\mathbb{N}^{d}}g_{J}^{(1)}(z)w^{J},\displaystyle\sum_{J\in\mathbb{N}^{d}}g_{J}^{(2)}(z)w^{J},\dots,\displaystyle\sum_{J\in\mathbb{N}^{d}}g_{J}^{\left(d\right)}(z)w^{J}\right),\\& f(z,w)=\displaystyle\sum_{J\in\mathbb{N}^{d}}f_{J}(z)w^{J}=\left(\displaystyle\sum_{J\in\mathbb{N}^{d}}f_{J}^{(1)}(z)w^{J},\displaystyle\sum_{J\in\mathbb{N}^{d}}f_{J}^{(2)}(z)w^{J},\dots,\displaystyle\sum_{J\in\mathbb{N}^{d}}f_{J}^{\left(N\right)}(z)w^{J}\right),\label{mapA} \end{split}
\end{equation}
where   $f_{J}(z)$ and $g_{J}(z)$ are formal (holomorphic) power series  without constant terms,  for all $J\in\mathbb{N}^{d}$.

In particular for $d=1$, we work with formal expansions 
\begin{equation}  g(z,w)= \displaystyle\sum_{k\in\mathbb{N}}g_{k}(z)w^{k},\quad\quad f(z,w)=\displaystyle\sum_{k\in\mathbb{N}}f_{k}(z)w^{k}=\left( \displaystyle\sum_{k\in\mathbb{N}}f_{k}^{(1)}(z)w^{k}, \displaystyle\sum_{k\in\mathbb{N}}f_{k}^{(2)}(z)w^{k},\dots, \displaystyle\sum_{k\in\mathbb{N}}f_{k}^{\left(N\right)}(z)w^{k}\right),\label{mapB}  
\end{equation}
where   $f_{k}(z)$ and $g_{k}(z)$ are formal power series without constant terms,  for all $k\in\mathbb{N}$. 
 
In order to understand the level of degeneracy of (\ref{mapA}) or (\ref{mapB}), we define
\begin{equation}\mbox{wt}\left\{z_{k}\right\}=1,\hspace{0.12 cm}\mbox{wt}\left\{\overline{z}_{k}\right\}=1,\hspace{0.12 cm} \mbox{wt}\left\{w_{l}\right\}=2,\quad \mbox{for all $k=1,\dots,N$ and $l=1,\dots,d$.}\label{weit}
\end{equation} 

The  weight of a formal power series in $(w,z)$  is defined to be the minimum of the weights of the  components from its formal expansion in weighted-homogeneous  polynomials. It is standard definition taken from \cite{V1} in order to purse further computations, assuming
\begin{equation}\begin{split}&\mbox{wt}\left\{\displaystyle\sum_{k\in\mathbb{N}}f_{k}^{\left(\alpha\right)}(z)w^{k}\right\}=p_{i},\quad\mbox{for all $\alpha\in\left[n_{i}, ,n_{i+1}-1\right)$, for all 
$i=1,\dots,N_{1}-1$,}\\& \mbox{wt}\left\{\displaystyle\sum_{k\in\mathbb{N}}g_{k}^{\left(\beta\right)}(z)w^{k}\right\}=q_{j},\quad\mbox{for all $\beta\in\left[m_{j},m_{j+1}-1\right)$, for all $j=1,\dots,d_{1}-1$,}\end{split}\label{teta1}
\end{equation} 
where we have used the following numbers 
\begin{equation}\begin{split}& \mathbb{N}^{\star} \supset \left\{p_{i}\right\}_{i=1,\dots,N_{1}},\hspace{0.1 cm} \left\{n_{i}\right\}_{i=1,\dots,N_{1}} ,\hspace{0.05 cm}    \quad\mbox{such that $n_{1}\leq n_{2}\leq \dots \leq  n_{N_{1}}$ with $n_{1}=1$ and $n_{N_{1}}=N$,} \\& \mathbb{N}^{\star} \supset \left\{q_{i}\right\}_{i=0,\dots,N'_{1}},\hspace{0.1 cm} \left\{m_{j}\right\}_{j=1,\dots,d_{1}},\quad   \mbox{such that $m_{1}\leq m_{2}\leq \dots \leq  m_{d_{1}}$ with  $m_{0}=1$ and $m_{d_{1}}=d$. }   \end{split}\label{teta2}
\end{equation} 

Since we work with   power series, we write
\begin{equation}\begin{split}& \hspace{0.17 cm} g_{J}^{\left(l\right)}(z)=\displaystyle\sum_{I=\left(i_{1},i_{2}, \dots,i_{N}\right)\in\mathbb{N}^{d}}a_{I}^{\left(l\right)}z_{1}^{i_{1}}z_{2}^{i_{2}}\dots z_{N-d}^{i_{N-d}}z_{N-d+1}^{i_{N-d+1}}\dots z_{N}^{i_{N}},\hspace{0.18 cm}\quad\mbox{for all $l= 1,\dots, d$ and $J=\left(j_{1},j_{2},\dots,j_{d}\right)\in\mathbb{N}^{d}$,}\\& f_{J}^{\left(l'\right)}(z)=\displaystyle\sum_{I=\left(i_{1},i_{2}, \dots,i_{N}\right)\in\mathbb{N}^{d}}b_{I}^{\left(l'\right)}z_{1}^{i_{1}}z_{2}^{i_{2}}\dots z_{N-d}^{i_{N-d}}z_{N-d+1}^{i_{N-d+1}}\dots z_{N}^{i_{N}},\quad\hspace{0.1 cm}\mbox{for all  $l'= 1,\dots, N$ and $J=\left(j_{1},j_{2},\dots,j_{d}\right)\in\mathbb{N}^{d}$.}\end{split} \label{ge1}
  \end{equation}

In particular for $d=1$, we write 
 \begin{equation}  g_{J}(z)=\displaystyle\sum_{I=\left(i_{1},i_{2}, \dots,i_{N}\right)\in\mathbb{N}^{d}}a_{I}z_{1}^{i_{1}}z_{2}^{i_{2}}\dots z_{N-d}^{i_{N-d}}z_{N-d+1}^{i_{N-d+1}}\dots z_{N}^{i_{N}},\quad \mbox{for all   $J=\left(j_{1},j_{2},\dots,j_{d}\right)\in\mathbb{N}^{d}$.} \label{ge11}
  \end{equation}

We move forward using the following: 
\subsection{Notations} If $F(z,\overline{z})$ is a formal power series in $(z,\overline{z})$, we define
\begin{equation*}  \left(\Re F\right)\left(z,\overline{z}\right)= \Re F\left(z,\overline{z}\right)=\frac{F\left(z,\overline{z}\right)+\overline{F\left(z,\overline{z}\right)}}{2},  \quad\quad   \left(\Im F\right)\left(z,\overline{z}\right)=\Im F\left(z,\overline{z}\right)=\frac{F\left(z,\overline{z}\right)-\overline{F\left(z,\overline{z}\right)}}{2\sqrt{-1}}.
\end{equation*} 

In particular, we have
\begin{equation}\begin{split}&  \hspace{0.1 cm}\left(\Re Q\right)(z,\overline{z})= \left(
\left(\Re Q_{1}\right)(z,\overline{z}) ,   \dots , \left(\Re Q_{d}\right)(z,\overline{z}) 
\right), \quad\hspace{0.25 cm}  \left(\Im Q\right)(z,\overline{z})= \left(
\left(\Im Q_{1}\right)(z,\overline{z}) ,  \dots , \left(\Im Q_{d}\right)(z,\overline{z}) 
\right),\\& \hspace{0.12 cm} \left(\Re \varphi\right)(z,\overline{z})=\left(
\left(\Re \varphi_{1}\right)(z,\overline{z})  , \dots , \left(\Re \varphi_{N}\right)(z,\overline{z})
\right),\quad\quad \left(\Im \varphi\right)(z,\overline{z})=\left(
\left(\Im \varphi_{1}\right)(z,\overline{z}) ,   \dots , \left(\Im \varphi_{N}\right)(z,\overline{z})
\right),\\&   \left(\Re\varphi'\right)(z,\overline{z})=\left(
\left(\Re\varphi'_{1}\right)(z,\overline{z}) ,   \dots , \left(\Re\varphi'_{N}\right)(z,\overline{z})
\right), \quad \hspace{0.04 cm} \left(\Im\varphi'\right)(z,\overline{z})=\left(\left(\Im\varphi'_{1}\right)(z,\overline{z}) ,   \dots , \left(\Im\varphi'_{N}\right)(z,\overline{z})
\right).\end{split}\label{Wla}
\end{equation} 

For any complex-valued smooth vector-function, denoted by $\rho\left(z,\overline{z},x\right):=\left(\rho_{1},\rho_{2},\dots,\rho_{d}\right)\left(z,\overline{z},x\right)$, 
we introduce the   matrix
\begin{equation} \gamma_{0}\left(z,\overline{z},x\right)=\begin{pmatrix}\frac{\partial \rho_{1}}{\partial z_{1}}\left(z,\overline{z},x\right)  & \frac{\partial \rho_{1}}{\partial z_{2}}\left(z,\overline{z},x\right) & \dots & \frac{\partial \rho_{1}}{\partial z_{N-d}}\left(z,\overline{z},x\right)\\ \frac{\partial \rho_{2}}{\partial z_{1}}\left(z,\overline{z},x\right)  & \frac{\partial \rho_{2}}{\partial z_{2}}\left(z,\overline{z},x\right) & \dots & \frac{\partial \rho_{2}}{\partial z_{N-d}}\left(z,\overline{z},x\right)\\ \vdots & \vdots & \ddots & \vdots     \\ \frac{\partial \rho_{d}}{\partial z_{1}}\left(z,\overline{z},x\right)  & \frac{\partial \rho_{d}}{\partial z_{2}}\left(z,\overline{z},x\right) & \dots & \frac{\partial \rho_{d}}{\partial z_{N-d}}\left(z,\overline{z},x\right)
\end{pmatrix} .\label{ro1}\end{equation}

In particular for $d=1$, we introduce the  matrix
\begin{equation} \tilde{\gamma}_{0}\left(z,\overline{z},x\right)=\begin{pmatrix}\frac{\partial \left(\Re \rho\right)}{\partial z_{1}}\left(z,\overline{z},x\right)  & \frac{\partial \left(\Re \rho\right)}{\partial z_{2}}\left(z,\overline{z},x\right) & \dots & \frac{\partial \left(\Re \rho\right)}{\partial z_{N-2}}\left(z,\overline{z},x\right)\\ \frac{\partial \left(\Im \rho\right)}{\partial z_{1}}\left(z,\overline{z},x\right)  & \frac{\partial \left(\Im \rho\right)}{\partial z_{2}}\left(z,\overline{z},x\right) & \dots & \frac{\partial \left(\Im \rho\right)}{\partial z_{N-2}}\left(z,\overline{z},x\right) 
\end{pmatrix} .\label{ro2}\end{equation}

Since we work with   power series, we write 
 \begin{equation}\begin{split}& \varphi_{l}(z,\overline{z})=\displaystyle\sum_{ I=\left(i_{1},i_{2}, \dots,i_{N}\right)\in\mathbb{N}^{d}\atop{I'=\left(i'_{1},i'_{2}, \dots,i'_{N}\right)\in\mathbb{N}^{d}}}\varphi^{(l)}_{I,I'}z_{1}^{i_{1}}z_{2}^{i_{2}}\dots z_{N-d}^{i_{N-d}}z_{N-d+1}^{i_{N-d+1}}\dots z_{N}^{i_{N}} \overline{z}_{1}^{i'_{1}}\overline{z}_{2}^{i'_{2}}\dots \overline{z}_{N-d}^{i'_{N-d}}\overline{z}_{N-d+1}^{i'_{N-d+1}}\dots \overline{z}_{N}^{i'_{N}} ,\quad\hspace{0.1 cm}\mbox{for all $l=1,\dots,d$,}\\&  \varphi'_{l}(z,\overline{z})=\displaystyle\sum_{ I=\left(i_{1},i_{2}, \dots,i_{N}\right)\in\mathbb{N}^{d}\atop{I'=\left(i'_{1},i'_{2}, \dots,i'_{N}\right)\in\mathbb{N}^{d}}}{\varphi'}^{(l)}_{I,I'}z_{1}^{i_{1}}z_{2}^{i_{2}}\dots z_{N-d}^{i_{N-d}}z_{N-d+1}^{i_{N-d+1}}\dots z_{N}^{i_{N}} \overline{z}_{1}^{i'_{1}}\overline{z}_{2}^{i'_{2}}\dots \overline{z}_{N-d}^{i'_{N-d}}\overline{z}_{N-d+1}^{i'_{N-d+1}}\dots \overline{z}_{N}^{i'_{N}} ,\quad\mbox{for all $l=1,\dots,d$.}\end{split}\label{ge2}
\end{equation} 

In particular for $d=1$, we write  
    \begin{equation}\begin{split}& \varphi(z,\overline{z})=\displaystyle\sum_{ I=\left(i_{1},i_{2}, \dots,i_{N}\right)\in\mathbb{N}^{d}\atop{I'=\left(i'_{1},i'_{2}, \dots,i'_{N}\right)\in\mathbb{N}^{d}}}\varphi_{I,I'}z_{1}^{i_{1}}z_{2}^{i_{2}}\dots z_{N-d}^{i_{N-d}}z_{N-d+1}^{i_{N-d+1}}\dots z_{N}^{i_{N}} \overline{z}_{1}^{i'_{1}}\overline{z}_{2}^{i'_{2}}\dots \overline{z}_{N-d}^{i'_{N-d}}\overline{z}_{N-d+1}^{i'_{N-d+1}}\dots \overline{z}_{N}^{i'_{N}},\\& \varphi'(z,\overline{z})=\displaystyle\sum_{ I=\left(i_{1},i_{2}, \dots,i_{N}\right)\in\mathbb{N}^{d}\atop{I'=\left(i'_{1},i'_{2}, \dots,i'_{N}\right)\in\mathbb{N}^{d}}}\varphi'_{I,I'}z_{1}^{i_{1}}z_{2}^{i_{2}}\dots z_{N-d}^{i_{N-d}}z_{N-d+1}^{i_{N-d+1}}\dots z_{N}^{i_{N}} \overline{z}_{1}^{i'_{1}}\overline{z}_{2}^{i'_{2}}\dots \overline{z}_{N-d}^{i'_{N-d}}\overline{z}_{N-d+1}^{i'_{N-d+1}}\dots \overline{z}_{N}^{i'_{N}} .\end{split}\label{ge21}
\end{equation} 

In particular, we write 
\begin{equation}z_{k}=\alpha_{k}+\sqrt{-1}\beta_{k},\quad\mbox{for all $k=1,\dots,N$}.
\label{alf}
\end{equation}

We move forward in order to study:

\section{CR Singularities in   Codimension at least $2$} 

Let $\delta>0$ be sufficiently small real number and 
$x=\left(x_{1},x_{2},\dots,x_{d}\right)$, for all  $x_{1},x_{2},\dots,x_{d}\in(0,\delta)$.  Restricting (\ref{M1}) when $w=x$, or equivalently when $w_{1}=x_{1}>0$, $w_{2}=x_{2}>0,\dots,w_{d}=x_{d}>0$,  we obtain by (\ref{Deg}), (\ref{real}) and (\ref{Degg})  a family of Real-Analytic Submanifolds     
 \begin{equation}M_{x}:\hspace{0.1 cm}\rho\left(z,\overline{z},x\right)=0,\quad\mbox{for}\hspace{0.1 cm}\rho\left(z,\overline{z},x\right):=   Q\left(z,\overline{z}\right)+\left(\Re\varphi\right)(z,\overline{z})-x. \label{98}\end{equation}
 
  Focusing on the real part in (\ref{M1}),   we assume $ \Im w= \left(\Im \varphi\right)(z,\overline{z}),
 $ or equivalently 
\begin{equation} \left(\left(\Im \varphi_{1}\right)(z,\overline{z})    ,\left(\Im \varphi_{2}\right)(z,\overline{z}),\dots,\left(\Im \varphi_{d}\right)(z,\overline{z})\right)=\left(\Im w_{1},\Im w_{2},\dots,\Im w_{d}\right).
 \label{asu}
\end{equation}

Now, let $z_{x}\in\mathbb{C}^{N}$ such that 
\begin{equation}\rho\left(z_{x},\overline{z_{x}},x\right)=0,\quad\quad\quad  \gamma_{0}\left(z_{x},\overline{z_{x}},x\right)\neq 0.\label{61}
\end{equation} 
 
From the Theorem of Implicit Functions applied in  (\ref{61}), it follows that it exists an analytic function, denoted by
\begin{equation} \psi:\mathbb{C}^{N}\times (0,\epsilon)\times(0,\epsilon)\dots \times(0,\epsilon) \rightarrow \mathbb{R}^{d},\quad  \psi:=\left(\psi_{N-d+1},\psi_{N-d+2},\dots,\psi_{N}\right),  \label{NX}
\end{equation}
 defined near $0\in\mathbb{C}^{N}$, where $\epsilon>0$ is chosen enough small, such that   
 \begin{equation}M_{x}:\hspace{0.1 cm}\Im \begin{pmatrix}z_{N-d+1}\\ z_{N-d+2}\\ \vdots
\\   z_{N}
\end{pmatrix}=\begin{pmatrix}\psi_{N-d+1}\left(x,z_{1},\overline{z}_{1},z_{2},\overline{z}_{2},\dots,z_{N-d},\overline{z}_{N-d},\Re z_{N-d+1},\Re z_{N-d+2},\dots, \Re z_{N}\right)\\ \psi_{N-d+2}\left(x,z_{1},\overline{z}_{1},z_{2},\overline{z}_{2},\dots,z_{N-d},\overline{z}_{N-d},\Re z_{N-d+1},\Re z_{N-d+2},\dots, \Re z_{N}\right)\\ \vdots
\\  \quad\quad \psi_{N}\left(x,z_{1},\overline{z}_{1},z_{2},\overline{z}_{2},\dots,z_{N-d},\overline{z}_{N-d},\Re z_{N-d+1},\Re z_{N-d+2},\dots, \Re z_{N}\right)
\end{pmatrix} .\label{NX1}\end{equation}  
 
We   rewrite the right-hand side from (\ref{NX1}) like 
\begin{equation} \psi_{k}=\psi_{0k}+\tilde{\psi}_{k},\quad\mbox{for all $k=1,\dots,d$,}\label{rez} 
\end{equation}
where we have used the following notations
\begin{equation}
\begin{split}& \psi_{0k}=\psi_{0k}\left(x,z_{1},z_{2},\dots,z_{N-d} \right) ,\quad \quad \quad \quad  \quad\quad \quad\quad \hspace{0.25 cm}  \quad \quad \quad\quad \quad \quad \quad \quad  \quad \quad \quad \quad    \quad\mbox{for all $k=1,\dots,d$,}\\&\hspace{0.1 cm} \tilde{\psi}_{k}= \tilde{\psi}\left(x,z_{1},\overline{z}_{1},z_{2},\overline{z}_{2},\dots,z_{N-d},\overline{z}_{N-d};\Re z_{N-d+1},\Re z_{N-d+2},\dots, \Re z_{N}\right),\quad   \quad\mbox{for all $k=1,\dots,d$,}\\& \hspace{0.1 cm} \psi_{k}= \tilde{\psi}\left(x,z_{1},\overline{z}_{1},z_{2},\overline{z}_{2},\dots,z_{N-d},\overline{z}_{N-d};\Re z_{N-d+1},\Re z_{N-d+2},\dots, \Re z_{N}\right),\quad     \quad\mbox{for all $k=1,\dots,d$,} \end{split}\label{rez1} 
\end{equation}
where we deal with:
\begin{itemize}
\item $\psi_{0k}$ is the analytic component in $\left(z_{1}, z_{2},\dots,z_{N-d} \right)$ and $x$ of $\psi_{k}$, for all $k=1,\dots,d$,     
\item $\tilde{\psi}_{k}$ is the analytic component
in   $x$,   $\left(z_{1}, z_{2},\dots,z_{N-d} \right)$, $\left(\overline{z}_{1}, \overline{z}_{2},\dots,\overline{z}_{N-d} \right)$ and $\left(\Re z_{N-d+1},\Re z_{N-d+2},\dots, \Re z_{N}\right)$ of $\psi_{k}$,      
for all $k=1,\dots,d$. \end{itemize}

We consider by (\ref{rez1}) the following change of coordinates
\begin{equation}\begin{split}& \quad\quad\quad\quad  \left(z'_{1}
, z'_{2}, \dots , z'_{N-d}
\right)=\left(z_{1}
, z_{2} , \dots , z_{N-d}
\right),\\&  \left(z'_{N-d+1}
,z'_{N-d+2},\dots,  z'_{N}
\right)=\left(z_{N-d+1},z_{N-d+2},  \dots,  z_{N}
\right)-\sqrt{-1}\left(\psi_{01},\psi_{02},\dots,\psi_{0d}\right)\left(x,  z_{1}
, z_{2} , \dots , z_{N-d}
  \right).  \end{split} \label{44}
\end{equation} 

In particular, (\ref{44})  implies
\begin{equation}\begin{split}& \quad\quad\quad\quad  \left(z_{1}
, z_{2}, \dots , z_{N-d}
\right)=\left(z'_{1}
, z'_{2} , \dots , z'_{N-d}
\right),\\&  \left(z_{N-d+1}
,z_{N-d+2},\dots,  z_{N}
\right)=\left(z'_{N-d+1},z'_{N-d+2},  \dots,  z'_{N}
\right)+\sqrt{-1}\left(\psi_{01},\psi_{02},\dots,\psi_{0d}\right)\left(x,  z'_{1}
, z'_{2} , \dots , z'_{N-d}
  \right).  \end{split} \label{441}
\end{equation}
 
 Then, (\ref{NX1}) becomes
 \begin{equation}  M_{x}:\hspace{0.1 cm}\Im \begin{pmatrix}z'_{N-d+1}\\ z'_{N-d+2}\\\vdots
  \\ z'_{N}
\end{pmatrix}=\begin{pmatrix}  A_{N-d+1}\left(x,z'_{1},z'_{2},\dots,z'_{N-d},\overline{z}'_{1},\overline{z}'_{2},\dots,\overline{z}'_{N-d}\right)\\ A_{N-d+2}\left(x,z'_{1},z'_{2},\dots,z'_{N-d},\overline{z}'_{1},\overline{z}'_{2},\dots,\overline{z}'_{N-d}\right)\\ \hspace{0.15 cm}\quad\quad\quad\quad\vdots
\\ \quad\quad\hspace{0.1 cm} A_{N}\left(x,z'_{1},z'_{2},\dots,z'_{N-d},\overline{z}'_{1},\overline{z}'_{2},\dots,\overline{z}'_{N-d}\right)
\end{pmatrix},\label{442}\end{equation}
where we have used the following notations
 \begin{equation}\begin{split} A_{N-d+k}= &\tilde{\psi}_{N-d+k}\left(x,z'_{1},z'_{2},\dots,z'_{N-d},\overline{z}'_{1},\overline{z}'_{2},\dots,\overline{z}'_{N-d},\right.\\& \quad\quad\quad\quad\hspace{0.1 cm} \left.\Re \left(z'_{N-d+1}+\psi_{N-d+1}\right),\dots,\Re \left(z'_{N}+\psi_{N}\right)\right)-\psi_{0k},\quad\mbox{for all $k=1,\dots,d$,}\end{split}\label{443} \end{equation} 
where we have used the following notations
 \begin{equation}
 \begin{split}& A_{N-d+k}=A_{N-d+k}\left(x,z'_{1},z'_{2},\dots,z'_{N-d},\overline{z}'_{1},\overline{z}'_{2},\dots,\overline{z}'_{N-d}\right),\quad\mbox{for all $k=1,\dots,d$,}\\& \hspace{0.04 cm} \psi_{N-d+k}=\psi_{N}\left(x, z'_{1}, z'_{2},\dots,z'_{N-d}\right),\quad\quad\quad\quad\quad\quad \quad\quad \quad\quad\quad\mbox{for all $k=1,\dots,d$,}\\& \quad\quad \hspace{0.07 cm}\psi_{0k}=\psi_{0k}\left(x,z'_{1},z'_{2},\dots,z'_{N-d}\right),\quad\quad\quad\quad\quad\quad\quad\quad\quad\quad\quad\mbox{for all $k=1,\dots,d$.}\end{split}\label{444}
 \end{equation}
  
In particular, we define \begin{equation}\mathbb{B}_{\epsilon}=\left\{z'=\left(z'_{1},z'_{2},\dots,z'_{N}\right)\in\mathbb{C}^{N}; \hspace{0.1 cm}z'_{1},z'_{2},\dots,z'_{N} \in \overline{\mathbb{B}\left(0,\epsilon\right)}\subset\mathbb{C}\right\},\label{233}
\end{equation}  
  
  In order to make  partial evaluations of convergence  of the formal equivalence (\ref{mapA}) in the  coordinates (\ref{441}), we use (\ref{NX})   respecting the first Segre mapping\cite{BREbook},\cite{Mir1},   defined by \begin{equation}v_{1}:\mathbb{C}^{N}\longrightarrow \mathbb{C}^{N},\quad\quad  v_{1}\left(z'_{1},z'_{2},\dots,z'_{N-d};z'_{N-d+1},\dots,z'_{N}\right)=\left(z'_{1},z'_{2},\dots,z'_{N-d};0,0,\dots,0\right).\label{23}\end{equation}

In particular for $d=1$, we have
\begin{equation}v_{1}:\mathbb{C}^{N}\longrightarrow \mathbb{C}^{N},\quad\quad  v_{1}\left(z'_{1},z'_{2},\dots,z'_{N-1}; z'_{N}\right)=\left(z'_{1},z'_{2},\dots,z'_{N-1};0\right).\label{231}\end{equation}

In order to understand   (\ref{M1}) by complexifying (\ref{442}), 
we   make evaluations using (\ref{23}) and the   family of Real-Analytic Submanifolds
$\left\{M_{x}\right\}_{x  \backsim 0\in\mathbb{R}^{d}}\subset \mathbb{C}^{N}$ defined as in (\ref{61}), (\ref{NX}) and (\ref{NX1}),  using by (\ref{23}) that  $v_{1}\left(z'\right) \in \mathcal{ M}_{x}$, where   the complexification of (\ref{442}) is denoted by $\left\{\mathcal{M}_{x}\right\}_{x  \backsim 0\in\mathbb{R}^{d}}\subset \mathbb{C}^{N}$. Then, (\ref{ge1})  implies
    \begin{equation}\begin{split}&  \hspace{0.18 cm}\left.g_{J}^{\left(l\right)}(z)\right|_{z'=v_{1}\left(z'\right)}=\left.\displaystyle\sum_{I=\left(i_{1},i_{2}, \dots,i_{N}\right)\in\mathbb{N}^{d}}a_{I}^{\left(l\right)}\left(\displaystyle\prod_{k=1}^{N-d}{z'}_{k}^{i_{k}}\right) \cdot \left(\displaystyle\prod_{k=1}^{d}\left(\sqrt{-1}A_{N-d+k}-\psi_{N-d+k} \right)^{i_{N-d+k}}\right)\right|_{z'=v_{1}\left(z'\right)},  \\& \left.f_{J}^{\left(l'\right)}(z)\right|_{z'=v_{1}\left(z'\right)}=\left.\displaystyle\sum_{I=\left(i_{1},i_{2}, \dots,i_{N}\right)\in\mathbb{N}^{d}}b_{I}^{\left(l'\right)}\left(\displaystyle\prod_{k=1}^{N-d}{z'}_{k}^{i_{k}}\right) \cdot \left(\displaystyle\prod_{k=1}^{d}\left(\sqrt{-1}A_{N-d+k}-\psi_{N-d+k} \right)^{i_{N-d+k}}\right)\right|_{z'=v_{1}\left(z'\right)}, \end{split} \label{ge111}
  \end{equation}
for all   $l=1,\dots,d$ and $l'=1,\dots,N$, such that (\ref{443}) and (\ref{444}) are satisfied. 

Next, (\ref{ge2})   implies
\begin{equation}\begin{split}  \left.\varphi_{l}\left(z,\overline{z} \right)\right|_{z'=v\left(z'\right)}=&\displaystyle\sum_{ I=\left(i_{1},i_{2}, \dots,i_{N}\right)\in\mathbb{N}^{d}\atop{I'=\left(i'_{1},i'_{2}, \dots,i'_{N}\right)\in\mathbb{N}^{d}}}\varphi^{(l)}_{I,I'} \left(\displaystyle\prod_{k=1}^{N-d}{z'}_{k}^{i_{k}}\right) \cdot \left(\displaystyle\prod_{k=1}^{d}\left(\sqrt{-1}A_{N-d+k}-\psi_{N-d+k} \right)^{i_{N-d+k}}\right)\cdot\\& \quad\quad\quad\quad\quad\quad\quad\quad\quad\quad\quad\quad  \left.\overline{\left(\displaystyle\prod_{k=1}^{N-d}{z'}_{k}^{i_{k}}\right) \cdot \left(\displaystyle\prod_{k=1}^{d}\left(\sqrt{-1}A_{N-d+k}-\psi_{N-d+k} \right)^{i_{N-d+k}}\right)}\right|_{z'=v_{1}\left(z'\right)} , \end{split}\label{ge31}
  \end{equation}
for all $l=1,\dots,d$, such that (\ref{443}) and (\ref{444}) are satisfied, and respectively 
 \begin{equation}\begin{split}&  \left.\varphi'_{l}\left(f(z,w),\overline{f(z,w)}\right)\right|_{z'=v_{1}\left(z'\right)}=\displaystyle\sum_{ I=\left(i_{1},i_{2}, \dots,i_{N}\right)\in\mathbb{N}^{d}\atop{I'=\left(i'_{1},i'_{2}, \dots,i'_{N}\right)\in\mathbb{N}^{d}}}{\varphi'}^{(l)}_{I,I'}\displaystyle\prod_{k=1}^{N-d}\left(\displaystyle\sum_{I=\left(i_{1},i_{2}, \dots,i_{N}\right)\in\mathbb{N}^{d}} b_{I}^{\left(k\right)}\left(\displaystyle\prod_{k=1}^{N-d}{z'}_{k}^{i_{k}}\right)\right. \cdot \\& \left.\quad\quad\quad\quad \left(\displaystyle\prod_{k=1}^{d}\left(\sqrt{-1}A_{N-d+k}-\psi_{N-d+k} \right)^{i_{N-d+k}}\right)\right)^{i_{k}}\cdot  \displaystyle\prod_{k=1}^{N-d}\left(\overline{\displaystyle\sum_{  I'=\left(i'_{1},i'_{2}, \dots,i'_{N}\right)\in\mathbb{N}^{d}}b_{I'}^{\left(k\right)}\left(\displaystyle\prod_{k=1}^{N-d}{z'}_{k}^{i'_{k}}\right)} \right. \cdot \\& \left. \quad\quad\quad\quad\quad\quad\quad\quad\quad\quad\quad\quad\quad\quad\quad\quad\quad\quad\quad\quad\quad\quad\quad\quad\left.\left(\overline{\displaystyle\prod_{k=1}^{d}\left(\sqrt{-1}A_{N-d+k}-\psi_{N-d+k} \right)^{i'_{N-d+k}}}\right)\right)^{i'_{k}}\right|_{z'=v_{1}\left(z'\right)}  ,\end{split} \label{ge32}
  \end{equation}
for all $l=1,\dots,d$, such that (\ref{443}) and (\ref{444}) are satisfied.
 
In the previous computations, we observe the existence of the expression  
  \begin{equation}Z_{J}(x,z,\overline{z})=\left( x_{1} +\sqrt{-1}\left( \Im \varphi_{1}\right)(z,\overline{z})\right)^{j_{1}}\dots \left( x_{d} +\sqrt{-1}\left( \Im \varphi_{d}\right)(z,\overline{z})\right)^{j_{d}},\quad\mbox{for all $J=\left(j_{1},j_{2}, \dots,j_{N}\right)\in\mathbb{N}^{d}$,}
  \label{zet}
  \end{equation}
such that (\ref{443}) and (\ref{444}) are satisfied.  
  
Now, we move forward using the   sums
 \begin{equation}\begin{split}&  S_{11}(x,z,\overline{z})=g_{J}^{\left(l\right)}(z)\left( x_{1} +\sqrt{-1}\left( \Im \varphi_{1}\right)(z,\overline{z})\right)^{j_{1}}\dots \left( x_{d} +\sqrt{-1}\left( \Im \varphi_{d}\right)(z,\overline{z})\right)^{j_{d}} , \\&  S_{12}(x,z,\overline{z})= \overline{g_{J}^{\left(l\right)}(z)}\left( x_{1} -\sqrt{-1}\left( \Im \varphi_{1}\right)(z,\overline{z})\right)^{j_{1}}\dots \left( x_{d} -\sqrt{-1}\left( \Im \varphi_{d}\right)(z,\overline{z})\right)^{j_{d}},  \\&  S_{21}(x,z,\overline{z})= f_{J}^{\left(l'\right)}(z)\left( x_{1} +\sqrt{-1}\left( \Im \varphi_{1}\right)(z,\overline{z})\right)^{j_{1}}\dots \left( x_{d} +\sqrt{-1}\left( \Im \varphi_{d}\right)(z,\overline{z})\right)^{j_{d}} ,  \\&  S_{22}(x,z,\overline{z})=\overline{f_{J}^{\left(l'\right)}(z)}\left( x_{1} -\sqrt{-1}\left( \Im \varphi_{1}\right)(z,\overline{z})\right)^{j_{1}}\dots \left( x_{d} -\sqrt{-1}\left( \Im \varphi_{d}\right)(z,\overline{z})\right)^{j_{d}} ,  \end{split}\label{sume1}
  \end{equation}  
for all   $l=1,\dots,d$ and $l'=1,\dots,N$, with disrespect to $J=\left(j_{1},j_{2}, \dots,j_{N}\right)\in\mathbb{N}^{d}$.   

It follows that
\begin{equation}\begin{split}& \left.S_{11}(x,z,\overline{z})\right|_{z'=v_{1}\left(z'\right)}= \left.g_{J}^{\left(l\right)}(z)\left( x_{1} +\sqrt{-1}\left( \Im \varphi_{1}\right)(z,\overline{z})\right)^{j_{1}}\dots \left( x_{d} +\sqrt{-1}\left( \Im \varphi_{d}\right)(z,\overline{z})\right)^{j_{d}}\right|_{z'=v_{1}\left(z'\right)}, \\& \left.S_{12}(x,z,\overline{z})\right|_{z'=v_{1}\left(z'\right)}=\left. \overline{g_{J}^{\left(l\right)}(z)}\left( x_{1} -\sqrt{-1}\left( \Im \varphi_{1}\right)(z,\overline{z})\right)^{j_{1}}\dots \left( x_{d} -\sqrt{-1}\left( \Im \varphi_{d}\right)(z,\overline{z})\right)^{j_{d}}\right|_{z'=v_{1}\left(z'\right)}
,  \\& \left.S_{21}(x,z,\overline{z})\right|_{z'=v_{1}\left(z'\right)}= \left.  f_{J}^{\left(l'\right)}(z)\left( x_{1} +\sqrt{-1}\left( \Im \varphi_{1}\right)(z,\overline{z})\right)^{j_{1}}\dots \left( x_{d} +\sqrt{-1}\left( \Im \varphi_{d}\right)(z,\overline{z})\right)^{j_{d}}\right|_{z'=v_{1}\left(z'\right)},  \\& \left.S_{22}(x,z,\overline{z})\right|_{z'=v_{1}\left(z'\right)}=  \left.\overline{f_{J}^{\left(l'\right)}(z)}\left( x_{1} -\sqrt{-1}\left( \Im \varphi_{1}\right)(z,\overline{z})\right)^{j_{1}}\dots \left( x_{d} -\sqrt{-1}\left( \Im \varphi_{d}\right)(z,\overline{z})\right)^{j_{d}}\right|_{z'=v_{1}\left(z'\right)},  \end{split}\label{sume1123}
  \end{equation}  
for all   $l=1,\dots,d$ and $l'=1,\dots,N$, with disrespect to $J=\left(j_{1},j_{2}, \dots,j_{N}\right)\in\mathbb{N}^{d}$.

It is desired to prove the following affirmation 
\begin{equation}\left.\frac{\partial^{I} g_{k}^{\left(l\right)}\left(z\right) }{\partial {z'}^{I} }\right\vert_{z'=v_{1}\left(z'\right)}\hspace{0.1 cm}\mbox{is convergent},\quad\mbox{for all $l=1,\dots,d$ and   $k\in\mathbb{N}$,  }\label{MN1}\end{equation}
for all $I\in\mathbb{N}^{N}$ with $\left|I\right|=n$   and $n\in\mathbb{N}$, and respectively  the following affirmation
\begin{equation}\left. \frac{\partial^{I} f^{\left(l'\right)}_{k}\left(z\right) }{\partial {z'}^{I} }\right\vert_{z'=v_{1}\left(z'\right)}\hspace{0.1 cm}\mbox{is convergent},\quad\mbox{for all  $l'=1,\dots,N$ and $k\in\mathbb{N}$,   }\label{MN2} 
\end{equation}
for all $I\in\mathbb{N}^{N}$ with $\left|I\right|=n$  and $n\in\mathbb{N}$. 

In particular for $d=1$, (\ref{MN1}) becomes
\begin{equation}\left.\frac{\partial^{I} g_{k}\left(z\right) }{\partial {z'}^{I} }\right\vert_{z'=v_{1}\left(z'\right)}\hspace{0.1 cm}\mbox{is convergent},\quad\mbox{for all  $k\in\mathbb{N}$,  }\label{MN1se}\end{equation}
for all $I\in\mathbb{N}^{N}$ with $\left|I\right|=n$     and $n\in\mathbb{N}$.

It remains to prove
\bl\label{tele1}   There exist   $\delta,\epsilon>0$ small enough and the  Real Submanifolds in $\mathbb{C}^{N}$ defined by
\begin{equation}\begin{split}& M_{x}:\hspace{0.1 cm}\mathbb{B}_{\epsilon}\cap\left\{z'=\left(z'_{1},z'_{2},\dots,z'_{N}\right)\in\mathbb{C}^{N};\hspace{0.1 cm} x=Q\left(z',\overline{z}'\right)+\left(\Re \varphi\right)\left(z',\overline{z}'\right)\right\}  ,\quad \mbox{for all $x=\left(x_{1},x_{2},\dots,x_{d}\right)\in\mathbb{R}_{+}^{d}$ with }\\&\quad\quad\quad\quad\quad\quad\quad\quad\quad\quad\quad\quad\quad\quad\quad\quad\quad\quad\quad\quad\quad\hspace{0.18 cm}\quad\quad\quad\quad\quad\quad\quad\quad\quad\quad\quad\quad\quad\quad\quad\quad  \mbox{$x_{1},x_{2},\dots,x_{d}$ near $\delta>0$.}\end{split}\label{rot}
\end{equation}

\el
\begin{proof}The Real Submanifolds (\ref{rot}) are defined by   (\ref{98}).   Then (\ref{rot}) has sense for $x_{1}=x_{2}=\dots=x_{d}=\delta>0$  enough small chosen such that (\ref{asu}) holds, because we can   eventually change the coordinates by multiplying with $-1$ in  (\ref{M1}). Furthermore, we can take values around $\delta>0$ for $x_{1},x_{2},\dots,x_{d}$ in order to assume that the set   (\ref{rot}) is not empty. It becomes clear by (\ref{alf}) the linear independence of  
\begin{equation*} \begin{split}&\left\{\left(\frac{\partial}{\partial \alpha_{1}}\left(Q_{u}\left(z',\overline{z}'\right)\right),\frac{\partial}{\partial \alpha_{2}}\left(Q_{u}\left(z',\overline{z}'\right)\right),\dots,\frac{\partial}{\partial \alpha_{N}}\left(Q_{u}\left(z',\overline{z}'\right)\right), \right.\right.\\&\left.\left.\quad \frac{\partial}{\partial \beta_{1}}\left(Q_{u}\left(z',\overline{z}'\right)\right),\frac{\partial}{\partial \beta_{2}}\left(Q_{u}\left(z',\overline{z}'\right)\right),\dots,\frac{\partial}{\partial \beta_{N}}\left(Q_{u}\left(z',\overline{z}'\right)\right)\right)\right\}_{u=1,\dots,d},\end{split}
\end{equation*}
because (\ref{Deg}) holds, implying the linear independence of  
\begin{equation*} \begin{split}&\left\{\left(\frac{\partial}{\partial \alpha_{1}}\left(\rho_{u}\left(z',\overline{z}'\right)\right),\frac{\partial}{\partial \alpha_{2}}\left(\rho_{u}\left(z',\overline{z}'\right)\right),\dots,\frac{\partial}{\partial \alpha_{N}}\left(\rho_{u}\left(z',\overline{z}'\right)\right), \right.\right.\\&\left.\left.\quad \frac{\partial}{\partial \beta_{1}}\left(\rho_{u}\left(z',\overline{z}'\right)\right),\frac{\partial}{\partial \beta_{2}}\left(\rho_{u}\left(z',\overline{z}'\right)\right),\dots,\frac{\partial}{\partial \beta_{N}}\left(\rho_{u}\left(z',\overline{z}'\right)\right)\right)\right\}_{u=1,\dots,d},\end{split}
\end{equation*}
by eventually taking $\epsilon>0$ small enough  and $\delta>0$ smaller, because 
 \begin{equation*}\begin{split} &\frac{\partial}{\partial \alpha_{k}}\left(\rho_{u}\left(z',\overline{z}'\right)\right)=\frac{\partial}{\partial \alpha_{k}}\left(Q_{u}\left(z',\overline{z}'\right)\right)+\mbox{O}\left(\epsilon\right),\quad  \mbox{for all $k=1,\dots,N$ and   $u=1,\dots,d$,} \\&
  \frac{\partial}{\partial \beta_{k}}\left(\rho_{u}\left(z',\overline{z}'\right)\right)=\frac{\partial}{\partial \beta_{k}}\left(Q_{u}\left(z',\overline{z}'\right)\right)+\mbox{O}\left(\epsilon\right),\quad  \mbox{for all $k=1,\dots,N$ and   $u=1,\dots,d$.}\end{split}   
\end{equation*}
\end{proof}

\bl   \label{tele11}  Let $x_{1},x_{2},\dots,x_{d}>0$ such that $f$ and $g$ are holomorphic near  the points of  
(\ref{98}). Then, $f$ and $g$ are holomorphic near $p=0\in\mathbb{C}^{N+d}$. 
\el  
 
 \begin{proof} The   defining equations  are considered  near   points $p_{x}\in M_{x}$. Then, we  write $f$ and $g$ as formal power series. Then,  $f$ and $g$ have  strictly positive radiuses of convergence near such points, or equivalently on the following compact set
 \begin{equation}M\cap \left\{(w,z)=\left(w;z_{1},z_{2},\dots,z_{N}\right)\in\mathbb{C}^{N};\hspace{0.1 cm}1-\delta<|w|^{2}+\left|z_{1}\right|^{2}+\left|z_{2}\right|^{2}+\dots+\left|z_{N}\right|^{2}<1+\delta\right\}, \label{op}
  \end{equation}
which can be covered by an infinite union of open sets respecting the considered  points $p_{x}\in M_{x}$.

Because this set (\ref{op}) is compact, it may be written as a finite union of such open sets, where  $f$ and $g$ are holomorphic functions.  The Phenomenon of Hartogs provides unique holomorphic extensions for these functions. We obtain the  convergence   near the origin.  
\end{proof}

We more forward to:  
 \section{CR Singularities in Codimension $2$: Case $Q\neq \overline{Q}$}
Let $\epsilon>$ be an enough small real number. We assume that   (\ref{real1}) holds and $d=1$.  Restricting (\ref{M2}) when   $w=x+\sqrt{-1}y$,  for  $x,y\in(0,\epsilon)$,  we obtain by (\ref{Deg}), (\ref{real1}) and (\ref{Degg1}) a family of Real-Analytic Submanifolds     
\begin{equation}N_{xy}:\hspace{0.1 cm}\rho\left(z,\overline{z},x,y\right)=0,\quad\mbox{for}\hspace{0.1 cm} \rho\left(z,\overline{z},x\right)= \left(\left(\Re Q\right)\left(z,\overline{z}\right)+\left(\Re\varphi\right)(z,\overline{z})-x, \left(\Im Q\right)\left(z,\overline{z}\right)+\left(\Im\varphi\right)(z,\overline{z})-y\right).\label{98q}\end{equation}

Now, let $z_{xy}\in\mathbb{C}^{N}$ such that 
\begin{equation}\rho\left(z_{xy},\overline{z_{xy}},x,y\right)=0,\quad\quad\quad  \tilde{\gamma}_{0}\left(z_{xy},\overline{z_{xy}},x,y\right)\neq 0.\label{61q}
\end{equation} 

From the Theorem of Implicit Functions applied in  (\ref{61}), it follows that there exists an analytic function,  denoted by
\begin{equation} \psi:\mathbb{C}^{N}\times (0,\epsilon)\times(0,\epsilon) \rightarrow \mathbb{R}^{2},\quad\quad\quad\quad\quad \psi:=\left(\psi_{N-1},\psi_{N}\right),  \label{NXq}
\end{equation}
for $\epsilon>0$ enough small chosen, defined near $0\in\mathbb{C}^{N}$, such that   
 \begin{equation}N_{xy}:\hspace{0.1 cm}\Im \begin{pmatrix}z_{N-1}
\\   z_{N}
\end{pmatrix}=\begin{pmatrix}\psi_{N-1}\left(x,y,z_{1},\overline{z}_{1},\dots,z_{N-2},\overline{z}_{N-2},\Re z_{N-1}, \Re z_{N}\right)
\\   \quad\psi_{N}\left(x,y,z_{1},\overline{z}_{1},\dots,z_{N-2},\overline{z}_{N-2},\Re z_{N-1}, \Re z_{N}\right)
\end{pmatrix}.\label{NX1q}\end{equation}  
 
We  then rewrite (\ref{NXq}) like
\begin{equation} \begin{pmatrix}\psi_{N-1}
\\   \quad\psi_{N}
\end{pmatrix} =\begin{pmatrix}\psi_{01}
\\   \psi_{02}
\end{pmatrix} + \begin{pmatrix}\tilde{\psi}_{1}
\\   \tilde{\psi}_{2}
\end{pmatrix},\label{rezq} 
\end{equation}   
where we have used the following notations
\begin{equation}
\begin{split}& \psi_{01}=\psi_{01}\left(x,z_{1},z_{2},\dots,z_{N-2} \right) ,\quad\quad 
\tilde{\psi}_{1}= \tilde{\psi}_{1}\left(x,z_{1},\overline{z}_{1},\dots,z_{N-2},\overline{z}_{N-2};\Re z_{N-1}, \Re z_{N}\right), \\& 
\psi_{02}=\psi_{02}\left(x,z_{1},z_{2},\dots,z_{N-2} \right) ,\quad\quad  
\tilde{\psi}_{2}= \tilde{\psi}_{2}\left(x,z_{1},\overline{z}_{1},\dots,z_{N-2},\overline{z}_{N-2};\Re z_{N-1}, \Re z_{N}\right).\end{split}\label{rez1q} 
\end{equation}
where we deal with:
\begin{itemize}
\item $\left(\psi_{01},\psi_{02}\right)$ is the analytic  component 
in $\left(z_{1}, z_{2},\dots,z_{N-2} \right)$, $x$ and $y$ of $\left(\psi_{1},\psi_{2},\right)$.
\item $\left(\tilde{\psi}_{1},\tilde{\psi}_{2}\right)$ is the analytic component 
in $x$, $y$, $\left(z_{1},\overline{z}_{1}, z_{2},\overline{z}_{2},\dots,z_{N-2},\overline{z}_{N-2}, \right)$ and $\left(\Re z_{N-1}, \Re z_{N}\right)$ of $\left(\psi_{1},\psi_{2},\right)$.  
\end{itemize}

We consider the following change of coordinates
\begin{equation}\begin{split}&   \left(z'_{1}
, z'_{2}, \dots , z'_{N-2}
\right)=\left(z_{1}
, z_{2} , \dots , z_{N-2}
\right),\\& \quad\quad\hspace{0.26 cm} \left(z'_{N-1}
,  z'_{N}
\right)=\left(z_{N-1},   z_{N}
\right)-\sqrt{-1}\left(\psi_{01},\psi_{02}\right)\left(x,y, z_{1}
, z_{2} , \dots , z_{N-2}
  \right).  \end{split} \label{44q}
\end{equation} 

In particular, (\ref{44q}) implies
\begin{equation} \begin{split}&\quad\hspace{0.1 cm}
z_{k}=z'_{k},\quad\quad\quad\quad\quad\quad\quad\quad\quad\quad\quad\quad\quad\quad\quad\quad\quad\quad\mbox{for all $k=1,\dots,N-2$},\\& z_{N-1}=z'_{N-1}+\sqrt{-1}\psi_{01}\left(x,y,z'_{1},z'_{2},\dots,z'_{N-2}\right), \\&\quad\hspace{0.01 cm} z_{N}=z'_{N}+\sqrt{-1}\psi_{02}\left(x,y,z'_{1},z'_{2},\dots,z'_{N-2}\right). \end{split} \label{441q}
\end{equation}

Then, (\ref{NXq}) becomes
 \begin{equation}  N'_{xy}:\hspace{0.1 cm}\Im \begin{pmatrix}z'_{N-1}
\\   z'_{N}
\end{pmatrix}=\begin{pmatrix} A\left(x,y,z'_{1}, z'_{2}, \dots , z'_{N-2}
 ,\overline{z}'_{1},z'_{2},   \dots , \overline{z}'_{N-2}\right)
\\ B\left(x,y,z'_{1}, z'_{2}, \dots , z'_{N-2}
 ,\overline{z}'_{1},z'_{2},   \dots , \overline{z}'_{N-2}\right)
\end{pmatrix}, \label{442q}\end{equation}
where we have used the following notations
 \begin{equation}\begin{split}&  \begin{pmatrix} A 
\\    B 
\end{pmatrix}:=\begin{pmatrix}\tilde{\psi}_{1}\left(x,y,z'_{1}, z'_{2}, \dots , z'_{N-2}
 ,\overline{z}'_{1},\overline{z}'_{2},   \dots , \overline{z}'_{N-2}
 ;  \Re \left(   z'_{N-1}+ \psi_{N-1} \right),  \Re \left(   z'_{N}
+ \psi_{N}\right) \right)-\Re\psi_{1}
\\   \tilde{\psi}_{2}\left(x,y,z'_{1}, z'_{2}, \dots , z'_{N-2}
 ,\overline{z}'_{1},\overline{z}'_{2},   \dots , \overline{z}'_{N-2}
 ;  \Re \left(   z'_{N-1}+ \psi_{N-1} \right),  \Re \left(   z'_{N}
+ \psi_{N}\right) \right)-\Re\psi_{2}
\end{pmatrix},\quad\mbox{for:}   \\& \hspace{0.1 cm}\quad\quad\quad\quad\quad\quad\quad\quad\quad   \psi_{1}=\psi_{1}\left(x,y, z'_{1},z'_{2} \dots,z'_{N-2}\right), \quad A=A\left(x,y,z'_{1}, z'_{2}, \dots , z'_{N-2}
 ,\overline{z}'_{1},z'_{2},   \dots , \overline{z}'_{N-2}\right),\\& \hspace{0.1 cm}\quad\quad\quad \quad\quad\quad\quad\quad\quad  \psi_{2}=\psi_{2}\left(x,y, z'_{1},z'_{2} \dots,z'_{N-2}\right),\quad   B=B\left(x,y,z'_{1}, z'_{2}, \dots , z'_{N-2}
 ,\overline{z}'_{1},z'_{2},   \dots , \overline{z}'_{N-2}\right).    \end{split} \label{443q}\end{equation}

Then, (\ref{NXq}) is  used in order to make  partial evaluations of convergence, of the formal equivalence (\ref{mapB}) in these coordinates (\ref{441}), respecting the first Segre mapping\cite{BREbook},\cite{Mir1},   defined by \begin{equation}\tilde{v}_{1}:\mathbb{C}^{N}\longrightarrow \mathbb{C}^{N},\quad\quad  \tilde{v}_{1}\left(z'_{1},z'_{2},\dots,z'_{N-2};z'_{N-1}, z'_{N}\right)=\left(z'_{1},z'_{2},\dots,z'_{N-2};0,0 \right).\label{23q}\end{equation}

Then, we   make evaluations using (\ref{23q}) and the   family of Real-Analytic Submanifolds
$ \left\{N_{xy}\right\}_{x,y  \backsim 0\in\mathbb{R}}\subset\mathbb{C}^{N}$ defined as in (\ref{61q}) and (\ref{NX1q}), in order to understand better (\ref{M2}) by complexifying (\ref{442q}), because we make evaluations by (\ref{23q}) previously the changing of coordinates by  (\ref{441q}) using by (\ref{23}) that $v_{1}\left(z'\right) \in \mathcal{ N}_{xy}$, 
where   the complexification of (\ref{442q}) is $
\left\{\mathcal{N}_{xy}\right\}_{x,y \backsim 0\in\mathbb{R}}\subset \mathbb{C}^{N}$  
Then,  (\ref{ge1})   implies
    \begin{equation}\begin{split}&  \left.\hspace{0.29 cm} g_{J}(z)\right|_{z'=\tilde{v}_{1}\left(z'\right)}=\left.\displaystyle\sum_{I=\left(i_{1},i_{2}, \dots,i_{N}\right)\in\mathbb{N}^{2}}a_{I}\left(\displaystyle\prod_{k=1}^{N-2}{z'}_{k}^{i_{k}}\right) \cdot \left( \left(   \sqrt{-1}A -\tilde{\psi}_{1}\right)^{i_{N-1}}\cdot \left(  \sqrt{-1}B -\tilde{\psi}_{2} \right)^{i_{N}}\right)\right|_{z'=\tilde{v}_{1}\left(z'\right)} ,  \\& \left.f_{J}^{\left(l'\right)}(z)\right|_{z'=\tilde{v}_{1}\left(z'\right)}=\left.\displaystyle\sum_{I=\left(i_{1},i_{2}, \dots,i_{N}\right)\in\mathbb{N}^{d}}b_{I}^{\left(l'\right)}\left(\displaystyle\prod_{k=1}^{N-2}{z'}_{k}^{i_{k}}\right) \cdot \left( \left(   \sqrt{-1}A -\tilde{\psi}_{1}\right)^{i_{N-1}}\cdot \left(  \sqrt{-1}B -\tilde{\psi}_{2} \right)^{i_{N}}\right)\right|_{z'=\tilde{v}_{1}\left(z'\right)},   \end{split} \label{ge111q}
  \end{equation}
for all $l'=1,\dots,N$, such that (\ref{443q}) holds, and respectively 
\begin{equation}\begin{split}\left. \varphi\left(z,\overline{z} \right)\right|_{z'=\tilde{v}\left(z'\right)}=&\displaystyle\sum_{ I=\left(i_{1},i_{2}, \dots,i_{N}\right)\in\mathbb{N}^{2}\atop{I'=\left(i'_{1},i'_{2}, \dots,i'_{N}\right)\in\mathbb{N}^{2}}}\varphi_{I,I'} \left(\displaystyle\prod_{k=1}^{N-2}{z'}_{k}^{i_{k}}\right) \cdot \left( \left(   \sqrt{-1}A -\tilde{\psi}_{1}\right)^{i_{N-1}}\cdot \left(  \sqrt{-1}B -\tilde{\psi}_{2} \right)^{i_{N}}\right)\cdot\\& \quad\quad\quad\quad\quad\quad\quad\quad\quad\quad\quad\quad  \left.\overline{\left(\displaystyle\prod_{k=1}^{N-2}{z'}_{k}^{i_{k}}\right) \cdot  \left(   \sqrt{-1}A -\tilde{\psi}_{1}\right)^{i_{N-1}}\cdot \left(  \sqrt{-1}B -\tilde{\psi}_{2} \right)^{i_{N}}}\right|_{z'=\tilde{v}\left(z'\right)}   , \end{split}\label{ge31q}
  \end{equation}
such that (\ref{443q}) holds, and respectively 
 \begin{equation}\begin{split}&  \left.\varphi'\left(f(z,w),\overline{f(z,w)}\right)\right|_{z'=\tilde{v}_{1}\left(z'\right)}=\displaystyle\sum_{ I=\left(i_{1},i_{2}, \dots,i_{N}\right)\in\mathbb{N}^{d}\atop{I'=\left(i'_{1},i'_{2}, \dots,i'_{N}\right)\in\mathbb{N}^{d}}}{\varphi'}_{I,I'}\displaystyle\prod_{k=1}^{N-2}\left(\displaystyle\sum_{I=\left(i_{1},i_{2}, \dots,i_{N}\right)\in\mathbb{N}^{d}} b_{I}^{\left(k\right)}\left(\displaystyle\prod_{k=1}^{N-2}{z'}_{k}^{i_{k}}\right)\right. \cdot \\& \left.\quad\quad\quad\quad\hspace{0.15 cm}\left(\left(   \sqrt{-1}A -\tilde{\psi}_{1}\right)^{i_{N-1}}\cdot \left(  \sqrt{-1}B -\tilde{\psi}_{2} \right)^{i_{N}}\right)\right)^{i_{k}}\cdot  \displaystyle\prod_{k=1}^{N-2}\left(\overline{\displaystyle\sum_{  I'=\left(i'_{1},i'_{2}, \dots,i'_{N}\right)\in\mathbb{N}^{d}}b_{I'}^{\left(k\right)}\left(\displaystyle\prod_{k=1}^{N-2}{z'}_{k}^{i'_{k}}\right)} \right. \cdot \\& \left. \quad\quad\quad\quad\quad\quad\quad\quad\quad\quad\quad\quad\quad\quad\quad\quad\quad\quad\quad\quad\quad\quad\quad\quad\quad\left. \overline{\left(   \sqrt{-1}A -\tilde{\psi}_{1}\right)^{i_{N-1}}\cdot \left(  \sqrt{-1}B -\tilde{\psi}_{2} \right)^{i_{N}}}\right)^{i'_{k}}\right|_{z'=\tilde{v}\left(z'\right)}  ,\end{split} \label{ge32q}
  \end{equation}
such that (\ref{443q}) holds.

In the previous computations, we observe the existence of the expression  
  \begin{equation}\tilde{Z}_{J}(x,y,z,\overline{z})=\left( x  +\sqrt{-1}\left( \Im \varphi\right)(z,\overline{z})\right)^{j_{1}}  \left( y +\sqrt{-1}\left( \Re \varphi\right)(z,\overline{z})\right)^{j_{2}},\quad\mbox{for all $J=\left(j_{1},j_{2}\right)\in\mathbb{N}^{2}$,}
  \label{zetX}
  \end{equation}
such that (\ref{442q}) and (\ref{443q}) are satisfied.  
  
Now, we move forward using the   sums
 \begin{equation}\begin{split}&  \tilde{S}_{11}(x,y,z,\overline{z})= g_{J}(z)\left( x  +\sqrt{-1}\left( \Im \varphi\right)(z,\overline{z})\right)^{j_{1}}  \left( y +\sqrt{-1}\left( \Re \varphi\right)(z,\overline{z})\right)^{j_{2}}, \\&  \tilde{S}_{12}(x,y,z,\overline{z})=  \overline{g_{J}(z)}\left( x  -\sqrt{-1}\left( \Im \varphi\right)(z,\overline{z})\right)^{j_{1}}  \left( y -\sqrt{-1}\left( \Re \varphi\right)(z,\overline{z})\right)^{j_{2}} ,  \\&  \tilde{S}_{21}(x,y,z,\overline{z})=  f_{J}^{\left(l'\right)}(z)\left( x  +\sqrt{-1}\left( \Im \varphi\right)(z,\overline{z})\right)^{j_{1}}  \left( y +\sqrt{-1}\left( \Re \varphi\right)(z,\overline{z})\right)^{j_{2}} ,  \\&  \tilde{S}_{22}(x,y,z,\overline{z})=\overline{f_{J}^{\left(l'\right)}(z)}\left( x  -\sqrt{-1}\left( \Im \varphi\right)(z,\overline{z})\right)^{j_{1}}  \left( y -\sqrt{-1}\left( \Re \varphi\right)(z,\overline{z})\right)^{j_{2}} ,  \end{split}\label{sume1X}
  \end{equation}  
for all   $l'=1,\dots,N$, with disrespect to $J=\left(j_{1},j_{2}\right)\in\mathbb{N}^{2}$.  

Now, we move forward using the   sums
 \begin{equation}\begin{split}& \left.\tilde{S}_{11}(x,y,z,\overline{z})\right|_{z'=v_{1}\left(z'\right)}= \left.g_{J}(z)\left( x  +\sqrt{-1}\left( \Im \varphi\right)(z,\overline{z})\right)^{j_{1}}  \left( y +\sqrt{-1}\left( \Re \varphi\right)(z,\overline{z})\right)^{j_{2}}\right|_{z'=v_{1}\left(z'\right)}, \\& \left.\tilde{S}_{12}(x,y,z,\overline{z})\right|_{z'=v_{1}\left(z'\right)}=\left. \overline{g_{J}(z)}\left( x  -\sqrt{-1}\left( \Im \varphi\right)(z,\overline{z})\right)^{j_{1}}  \left( y -\sqrt{-1}\left( \Re \varphi\right)(z,\overline{z})\right)^{j_{2}}\right|_{z'=v_{1}\left(z'\right)}
,  \\& \left.\tilde{S}_{21}(x,y,z,\overline{z})\right|_{z'=v_{1}\left(z'\right)}= \left.  f_{J}^{\left(l'\right)}(z)\left( x  +\sqrt{-1}\left( \Im \varphi\right)(z,\overline{z})\right)^{j_{1}}  \left( y +\sqrt{-1}\left( \Re \varphi\right)(z,\overline{z})\right)^{j_{2}}\right|_{z'=v_{1}\left(z'\right)},  \\& \left.\tilde{S}_{22}(x,y,z,\overline{z})\right|_{z'=v_{1}\left(z'\right)}=  \left.\overline{f_{J}^{\left(l'\right)}(z)}\left( x  -\sqrt{-1}\left( \Im \varphi\right)(z,\overline{z})\right)^{j_{1}}  \left( y -\sqrt{-1}\left( \Re \varphi\right)(z,\overline{z})\right)^{j_{2}}\right|_{z'=v_{1}\left(z'\right)},  \end{split}\label{sume1X123}
  \end{equation}  
for all   $l'=1,\dots,N$, with disrespect to $J=\left(j_{1},j_{2}\right)\in\mathbb{N}^{2}$.

  It remains to prove
\bl\label{tele2}   There exist   $\delta,\epsilon>0$ small enough and the  Real Submanifolds  in $\mathbb{C}^{N}$ defined by
\begin{equation}  N_{xy}:\hspace{0.1 cm}\mathbb{B}_{\epsilon}\cap\left\{z'=\left(z'_{1},z'_{2},\dots,z'_{N}\right)\in\mathbb{C}^{N};\quad x=\left(\Re Q\right)\left(z',\overline{z}'\right)+\left(\Re \varphi\right)\left(z',\overline{z}'\right)\atop{\hspace{0.14 cm}\quad\quad\quad\quad\quad\quad\quad\quad\quad\quad\quad\hspace{0.35 cm} y=\left(\Im Q\right)\left(z',\overline{z}'\right)+\left(\Im \varphi\right)\left(z',\overline{z}'\right)}\right\}  ,\quad\mbox{for all $x,y\in\mathbb{R}_{+}$ with  $x,y$ near $\delta>0$.}\label{rotse}
\end{equation}

\el
\begin{proof}The Real Submanifolds (\ref{rotse}) are defined by  (\ref{98}) near (\ref{pee1}). Then (\ref{rotse}) has sense for $x=y=\delta>0$ small enough taken, because we can   eventually change the coordinates by a simple rotation in   (\ref{M2}). Furthermore, we can take values around $\delta>0$ for $x,y$ in order to assume that  (\ref{rotse}) is not empty. Then, (\ref{Deg})   implies by (\ref{alf}) the linear independence of 
\begin{equation*}\left\{\begin{split}&  \left(\frac{\partial}{\partial \alpha_{1}}\left(\rho_{1}\left(z',\overline{z}'\right)\right),\frac{\partial}{\partial \alpha_{2}}\left(\rho_{1}\left(z',\overline{z}'\right)\right),\dots,\frac{\partial}{\partial \alpha_{N}}\left(\rho_{1}\left(z',\overline{z}'\right)\right),\right.\\&\left. \quad \frac{\partial}{\partial \beta_{1}}\left(\rho_{1}\left(z',\overline{z}'\right)\right),\frac{\partial}{\partial \beta_{2}}\left(\rho_{1}\left(z',\overline{z}'\right)\right),\dots,\frac{\partial}{\partial \beta_{N}}\left(\rho_{1}\left(z',\overline{z}'\right)\right)\right),  \\&   \left(\frac{\partial}{\partial \alpha_{1}}\left(\rho_{2}\left(z',\overline{z}'\right)\right),\frac{\partial}{\partial \alpha_{2}}\left(\rho_{2}\left(z',\overline{z}'\right)\right),\dots,\frac{\partial}{\partial \alpha_{N}}\left(\rho_{2}\left(z',\overline{z}'\right)\right), \right.\\&\quad \left. \frac{\partial}{\partial \beta_{1}}\left(\rho_{2}\left(z',\overline{z}'\right)\right),\frac{\partial}{\partial \beta_{2}}\left(\rho_{2}\left(z',\overline{z}'\right)\right),\dots,\frac{\partial}{\partial \beta_{N}}\left(\rho_{2}\left(z',\overline{z}'\right)\right)\right)     \end{split}\right\},
\end{equation*}
by eventually taking $\epsilon>0$ small enough  and $\delta>0$ smaller, because 
\begin{equation*}\begin{split}&\frac{\partial}{\partial \alpha_{k}}\left(\rho_{1}\left(z',\overline{z}'\right)\right)=\frac{\partial}{\partial \alpha_{k}}\left(\left(\Re Q\right)\left(z',\overline{z}'\right)\right)+\mbox{O}\left(\epsilon\right),\quad  \mbox{for all $k=1,\dots,N$,}\\& \frac{\partial}{\partial \beta_{k}}\left(\rho_{1}\left(z',\overline{z}'\right)\right)=\frac{\partial}{\partial \beta_{k}}\left(\left(\Re Q\right)\left(z',\overline{z}'\right)\right)+\mbox{O}\left(\epsilon\right),\quad  \mbox{for all $k=1,\dots,N$,}\\& \frac{\partial}{\partial \alpha_{k}}\left(\rho_{2}\left(z',\overline{z}'\right)\right)=\frac{\partial}{\partial \alpha_{k}}\left(\left(\Im Q\right)\left(z',\overline{z}'\right)\right)+\mbox{O}\left(\epsilon\right),\quad \mbox{for all $k=1,\dots,N$,} \\&    \frac{\partial}{\partial \beta_{k}}\left(\rho_{2}\left(z',\overline{z}'\right)\right)=\frac{\partial}{\partial \beta_{k}}\left(\left(\Im Q\right)\left(z',\overline{z}'\right)\right)+\mbox{O}\left(\epsilon\right),\quad  \mbox{for all $k=1,\dots,N$.}\end{split}
\end{equation*}
\end{proof}

\bl\label{tele21}  Let $x,y>0$ such that $f$ and $g$ are holomorphic near the points of 
(\ref{98q}). Then, $f$ and $g$ are holomorphic near $p=0\in\mathbb{C}^{N+1}$.
\el  
 
 \begin{proof}
The   defining equations  are considered  near  points $p_{xy}\in N_{xy}$. Then, we  write $f$ and $g$ as formal power series. Then,  $f$ and $g$ have  strictly positive radiuses of convergence near such points, or equivalently on the following compact set
 \begin{equation}M\cap \left\{(w,z)=\left(w;z_{1},z_{2},\dots,z_{N}\right)\in\mathbb{C}^{N};\hspace{0.1 cm}1-\delta<|w|^{2}+\left|z_{1}\right|^{2}+\left|z_{2}\right|^{2}+\dots+\left|z_{N}\right|^{2}<1+\delta\right\}, \label{p}
  \end{equation}
which can be covered by an infinite union of open sets respecting the considered  points $p_{xy}\in N_{xy}$.

Because this set (\ref{p}) is compact, it may be written as a finite union of such open sets, where  $f$ and $g$ are holomorphic functions.  The Phenomenon of Hartogs provides unique holomorphic extensions for these functions. Their  convergence becomes clear near the origin. 
 \end{proof}

Going forward,  we construct:
\section{General Analytic Systems }

 Let $M,M'\subset\mathbb{C}^{N+1}$ be two  Real-Analytic Submanifolds   as in (\ref{W}) such that it exists a formal mapping  like (\ref{mapA}) satisfying the conditions of Theorem \ref{t1}  in respect to (\ref{coord}).  We obtain
\begin{equation}\begin{split}  \displaystyle\sum_{ J=\left(j_{1},j_{2}, \dots,j_{d}\right)\in \mathbb{N}^{d}\atop{\left|J\right| \geq 1}}g_{J}^{\left(l\right)}(z)w_{1}^{j_{1}}w_{2}^{j_{2}} \dots w_{d}^{j_{d}}=&\left(Q'_{l}+\varphi'_{l}\right)\left( \displaystyle\sum_{ J=\left(j_{1},j_{2}, \dots,j_{d}\right)\in \mathbb{N}^{d}\atop{\left|J\right| \geq 1}}f_{J}^{\left(l\right)}w_{1}^{j_{1}}w_{2}^{j_{2}} \dots w_{d}^{j_{d}},\right.\\&  \quad\quad\quad\quad\left.\overline{ \displaystyle\sum_{ J=\left(j_{1},j_{2}, \dots,j_{d}\right)\in \mathbb{N}^{d}\atop{\left|J\right| \geq 1}}f_{J}(z)w_{1}^{j_{1}}w_{2}^{j_{2}} \dots w_{d}^{j_{d}}}\right),\quad\mbox{for all $l=1,\dots,d$.} \end{split} \label{V1}\end{equation}

Separating the real parts from the imaginary parts in (\ref{V1}), we obtain
\begin{equation}\begin{split}&  \frac{ \displaystyle\sum_{ J=\left(j_{1},j_{2}, \dots,j_{d}\right)\in \mathbb{N}^{d}\atop{\left|J\right| \geq 1}}\left(g_{J}^{\left(l\right)}(z)Z_{J}(x,z,\overline{z})-\overline{g_{J}^{\left(l\right)}(z)Z_{J}(x,z,\overline{z})}\right)}{2\sqrt{-1}} =  \left(\Im\varphi'_{l}\right)\left(\displaystyle\sum_{ J=\left(j_{1},j_{2}, \dots,j_{d}\right)\in \mathbb{N}^{d}\atop{\left|J\right| \geq 1}}f_{J}(z)Z_{J}(x,z,\overline{z}),\right.\\& \quad\quad\quad\quad\quad\quad\quad\quad\quad\quad\quad\quad\quad\quad\quad\quad\quad\quad\quad\quad\quad\quad\quad\quad\quad\quad\quad\quad\quad\quad\quad\hspace{0.18 cm} \left. \overline{ \displaystyle\sum_{ J=\left(j_{1},j_{2}, \dots,j_{d}\right)\in \mathbb{N}^{d}\atop{\left|J\right| \geq 1}}f_{J}(z)Z_{J}(x,z,\overline{z})} \right),\quad\mbox{for all $l=1,\dots,d$,} \end{split} \label{V21}\end{equation}
and respectively, we obtain
\begin{equation}\begin{split}&  \frac{\displaystyle\sum_{ J=\left(j_{1},j_{2}, \dots,j_{d}\right)\in \mathbb{N}^{d}\atop{\left|J\right| \geq 1}}\left(g_{J}^{\left(l\right)}(z)Z_{J}(x,z,\overline{z})-\overline{g_{J}^{\left(l\right)}(z)Z_{J}(x,z,\overline{z})}\right)}{2}= \left(Q'_{l}+\Im\varphi'_{l}\right)\left(\displaystyle\sum_{ J=\left(j_{1},j_{2}, \dots,j_{d}\right)\in \mathbb{N}^{d}\atop{\left|J\right| \geq 1}}f_{J}(z)Z_{J}(x,z,\overline{z}),\right.\\& \quad\quad\quad\quad\quad\quad\quad\quad\quad\quad\quad\quad\quad\quad\quad\quad\quad\quad\quad\quad\quad\quad\quad\quad\quad\quad\quad\quad\quad\quad\quad\quad\hspace{0.23 cm}  \left. \overline{ \displaystyle\sum_{ J=\left(j_{1},j_{2}, \dots,j_{d}\right)\in \mathbb{N}^{d}\atop{\left|J\right| \geq 1}}f_{J}(z)Z_{J}(x,z,\overline{z})} \right),\quad\mbox{for all $l=1,\dots,d$,} \end{split} \label{V22}\end{equation}
 when (\ref{98}) and (\ref{asu}) hold.   

The equations   (\ref{V21}) and  (\ref{V22}) are used   for suitable $\delta,\epsilon>0$ provided by Lemma \ref{tele11}. In particular, we chose by (\ref{98}) and (\ref{asu}) a point $z_{x}\in\mathbb{C}^{N}$ satisfying (\ref{61}).  Then, we    use  the Real-Analytic Submanifold $M_{x}$, which are resulted by (\ref{44}) and (\ref{441}) from The Implicit Function Theorem. We observe that \begin{equation}\overline{\frac{\partial g}{\partial z_{l}}(0)}\hspace{0.1 cm} \mbox{and} \hspace{0.1 cm}  \overline{\frac{\partial f_{k}}{\partial z_{l}}(0)} \hspace{0.1 cm}  \mbox{may not vanish, for all $k\in\mathbb{N}$ and $l=1,\dots,N$}. \label{V5}\end{equation}

Next, we identify the coefficients of $x$ in   (\ref{V21}) and (\ref{V22}).  We obtain that (\ref{MN2}) holds for $n=0$ and   the following analytic system:
\begin{equation} \begin{split}&\hspace{0.18 cm}  f_{1}\left(v_{1}\left(z'\right),x\right)+  P_{1}\left(f_{1}\left(v_{1}\left(z'\right),x\right),\dots,f_{N}\left(v_{1}\left(z'\right),x\right)\right)  =\chi_{1}\left(v_{1}\left(z'\right),x\right),\\&\hspace{0.18 cm} f_{2}\left(v_{1}\left(z'\right),x\right)   + P_{2}\left(f_{1}\left(v_{1}\left(z'\right),x\right),\dots,f_{N}\left(v_{1}\left(z'\right),x\right)\right) =\chi_{2}\left(v_{1}\left(z'\right),x\right), \\& \dots\dots\dots\dots\dots\dots\dots\dots\dots\dots\dots\dots\dots\dots\dots\dots\dots\dots\dots\dots\dots\dots\dots\dots\dots    \\&    f_{N}\left(v_{1}\left(z'\right),x\right)  +P_{N}\left(f_{1}\left(v_{1}\left(z'\right),x\right),\dots,f_{N}\left(v_{1}\left(z'\right),x\right)\right) =\chi_{N}\left(v_{1}\left(z'\right),x\right).  \end{split} \label{V6}\end{equation}  
  
 We apply Proposition $4.2$ from Mir\cite{Mir1} in (\ref{V6}). Then, (\ref{MN1}) holds for $n=0$.

Next, we can take derivatives with respect to $z_{1},z_{2},\dots,z_{N-d}$ in (\ref{V22}) of length $\left|I\right|=n$, for $I\in\mathbb{N}^{N}$  and a fixed $n\in\mathbb{N}$.  Repeating the previous procedure evaluating the resulted identity using (\ref{NX}), we conclude  the following analytic system
\begin{equation} \begin{split}&\hspace{0.17 cm} \frac{\partial^{I} f_{1}}{\partial z^{I}}\left(v_{1}\left(z'\right),x\right) + P'_{1}\left(\frac{\partial^{I} f_{1}}{\partial z^{I}}\left(v_{1}\left(z'\right),x\right),\dots,\frac{\partial^{I} f_{N}}{\partial z^{I}}\left(v_{1}\left(z'\right),x\right)\right) =\chi'_{1}\left(v_{1}\left(z'\right),x\right),  \\&\quad\hspace{0.15 cm} \frac{\partial^{I} f_{2}}{\partial z^{I}}\left(v_{1}\left(z'\right),x\right)  + P'_{2}\left(\frac{\partial^{I} f_{1}}{\partial z^{I}}\left(v_{1}(z),x\right),\dots,\frac{\partial^{I} f_{N}}{\partial z^{I}}\left(v_{1}\left(z'\right),x\right)\right)  =\chi'_{2}\left(v_{1}(z),x\right),   \\&\dots\dots\dots\dots\dots\dots\dots\dots\dots\dots\dots\dots\dots\dots\dots\dots\dots\dots\dots\dots\dots\dots\dots\dots\dots\dots\dots\dots\dots \\&\frac{\partial^{I}f_{N}}{\partial z^{I}}\left(v_{1}\left(z'\right),x\right)  + P'_{N}\left(\frac{\partial^{I} f_{1}}{\partial z^{I}}\left(v_{1}\left(z'\right),x\right),\dots,\frac{\partial^{I} f_{N}}{\partial z^{I}}\left(v_{1}\left(z'\right),x\right)\right) =\chi_{N}'\left(v_{1}\left(z'\right),x\right),  \end{split} \label{222CC}\end{equation}

 We apply Proposition $4.2$ from Mir\cite{Mir1} in (\ref{222CC}).   Then, (\ref{MN1}) holds for $n\in\mathbb{N}$ respecting an induction procedure.
 
\section{Analytic Systems in Codimension  $4$}
 
  Let $M,M'\subset\mathbb{C}^{N+1}$ be two  Real-Analytic Submanifolds   as in  (\ref{W}) such that it exists a formal mapping  like (\ref{mapA}) satisfying the conditions of Theorem \ref{t6}  in respect to (\ref{coord}).  We obtain
\begin{equation}   \displaystyle\sum_{ J=\left(j_{1},j_{2}\right)\in \mathbb{N}^{2}\atop{\left|J\right| \geq 1}}g_{J}^{\left(l\right)}(z)w_{1}^{j_{1}}w_{2}^{j_{2}} = \left(Q'_{l}+\varphi'_{l}\right)\left( \displaystyle\sum_{ J=\left(j_{1},j_{2}\right)\in \mathbb{N}^{2}\atop{\left|J\right| \geq 1}}f_{J}^{\left(l\right)}w_{1}^{j_{1}}w_{2}^{j_{2}} , \overline{ \displaystyle\sum_{ J=\left(j_{1},j_{2}\right)\in \mathbb{N}^{d}\atop{\left|J\right| \geq 1}}f_{J}(z)w_{1}^{j_{1}}w_{2}^{j_{2}} }\right),\quad\mbox{for all $l=1,2$.}   \label{V1bu1}\end{equation}

Separating the real parts from the imaginary parts in (\ref{V1bu1}), we obtain
\begin{equation}\begin{split}&  \frac{ \displaystyle\sum_{ J=\left(j_{1},j_{2}\right)\in \mathbb{N}^{2}\atop{\left|J\right| \geq 1}}\left(g_{J}^{\left(l\right)}(z)Z_{J}(x,y,z,\overline{z})-\overline{g_{J}^{\left(l\right)}(z)Z_{J}(x,y,z,\overline{z})}\right)}{2\sqrt{-1}} =  \left(\Im\varphi'_{l}\right)\left(\displaystyle\sum_{ J=\left(j_{1},j_{2}\right)\in \mathbb{N}^{2}\atop{\left|J\right| \geq 1}}f_{J}(z)Z_{J}(x,y,z,\overline{z}),\right.\\& \quad\quad\quad\quad\quad\quad\quad\quad\quad\quad\quad\quad\quad\quad\quad \quad\quad\quad\quad\quad\quad\quad\quad\quad\quad\quad\quad\quad\quad\quad\hspace{0.18 cm} \left. \overline{ \displaystyle\sum_{ J=\left(j_{1},j_{2}\right)\in \mathbb{N}^{2}\atop{\left|J\right| \geq 1}}f_{J}(z)Z_{J}(x,y,z,\overline{z})} \right),\quad\mbox{for all $l=1,2$,} \end{split} \label{V21bu1}\end{equation}
and respectively, we obtain
\begin{equation}\begin{split}&  \frac{\displaystyle\sum_{ J=\left(j_{1},j_{2}\right)\in \mathbb{N}^{2}\atop{\left|J\right| \geq 1}}\left(g_{J}^{\left(l\right)}(z)Z_{J}(x,y,z,\overline{z})-\overline{g_{J}^{\left(l\right)}(z)Z_{J}(x,y,z,\overline{z})}\right)}{2}= \left(Q'_{l}+\Im\varphi'_{l}\right)\left(\displaystyle\sum_{ J=\left(j_{1},j_{2}\right)\in \mathbb{N}^{2}\atop{\left|J\right| \geq 1}}f_{J}(z)Z_{J}(x,y,z,\overline{z}),\right.\\& \quad\quad\quad\quad\quad\quad\quad\quad\quad\quad\quad\quad\quad\quad\quad\quad\quad\quad\quad\quad\quad\quad\quad\quad\quad\quad\quad\quad\quad\quad\quad\quad\hspace{0.23 cm}  \left. \overline{ \displaystyle\sum_{ J=\left(j_{1},j_{2}\right)\in \mathbb{N}^{2}\atop{\left|J\right| \geq 1}}f_{J}(z)Z_{J}(x,y,z,\overline{z})} \right),\quad\mbox{for all $l=1,2$,} \end{split} \label{V22bu1}\end{equation}
 when (\ref{98q}) holds.   

The equations   (\ref{V21bu1}) and  (\ref{V22bu1}) are used   for suitable $\delta,\epsilon>0$ provided by Lemma \ref{tele21}. In particular, we chose by (\ref{98q}) a point $z_{x}\in\mathbb{C}^{N}$ satisfying (\ref{61q}).  Then, we    use  the Real-Analytic Submanifold $M_{x}$, which are resulted by (\ref{44q}) and (\ref{441q}) from The Implicit Function Theorem.

Next, we identify the coefficients of $x$ in   (\ref{V21bu1}) and  (\ref{V22bu1}).  We obtain that (\ref{MN2}) holds for $n=0$ and   the following analytic system:
\begin{equation} \begin{split}&\hspace{0.18 cm}  f_{1}\left(v_{1}\left(z'\right),x,y\right)+  P_{1}\left(f_{1}\left(v_{1}\left(z'\right),x,y\right),\dots,f_{N}\left(v_{1}\left(z'\right),x,y\right)\right)  =\chi_{1}\left(v_{1}\left(z'\right),x,y\right),\\&\hspace{0.18 cm} f_{2}\left(v_{1}\left(z'\right),x,y\right)   + P_{2}\left(f_{1}\left(v_{1}\left(z'\right),x,y\right),\dots,f_{N}\left(v_{1}\left(z'\right),x,y\right)\right) =\chi_{2}\left(v_{1}\left(z'\right),x,y\right), \\& \dots\dots\dots\dots\dots\dots\dots\dots\dots\dots\dots\dots\dots\dots\dots\dots\dots\dots\dots\dots\dots\dots\dots\dots\dots    \\&    f_{N}\left(v_{1}\left(z'\right),x,y\right)  +P_{N}\left(f_{1}\left(v_{1}\left(z'\right),x,y\right),\dots,f_{N}\left(v_{1}\left(z'\right),x,y\right)\right) =\chi_{N}\left(v_{1}\left(z'\right),x,y\right).  \end{split} \label{V6bu1}\end{equation}  
  
 We apply Proposition $4.2$ from Mir\cite{Mir1} in (\ref{V6}). Then, (\ref{MN1}) holds for $n=0$.

Next, we can take derivatives with respect to $z_{1},z_{2},\dots,z_{N-2}$ in (\ref{V22}) of length $\left|I\right|=n$, for $I\in\mathbb{N}^{N}$  and a fixed $n\in\mathbb{N}$.  Repeating the previous procedure evaluating the resulted identity using (\ref{NX}), we conclude  the following analytic system
\begin{equation} \begin{split}&\hspace{0.17 cm} \frac{\partial^{I} f_{1}}{\partial z^{I}}\left(v_{1}\left(z'\right),x,y\right) + P'_{1}\left(\frac{\partial^{I} f_{1}}{\partial z^{I}}\left(v_{1}\left(z'\right),x,y\right),\dots,\frac{\partial^{I} f_{N}}{\partial z^{I}}\left(v_{1}\left(z'\right),x,y\right)\right) =\chi'_{1}\left(v_{1}\left(z'\right),x,y\right),  \\&\quad\hspace{0.15 cm} \frac{\partial^{I} f_{2}}{\partial z^{I}}\left(v_{1}\left(z'\right),x,y\right)  + P'_{2}\left(\frac{\partial^{I} f_{1}}{\partial z^{I}}\left(v_{1}(z),x,y\right),\dots,\frac{\partial^{I} f_{N}}{\partial z^{I}}\left(v_{1}\left(z'\right),x,y\right)\right)  =\chi'_{2}\left(v_{1}(z),x,y\right),   \\&\dots\dots\dots\dots\dots\dots\dots\dots\dots\dots\dots\dots\dots\dots\dots\dots\dots\dots\dots\dots\dots\dots\dots\dots\dots\dots\dots\dots\dots \\&\frac{\partial^{I}f_{N}}{\partial z^{I}}\left(v_{1}\left(z'\right),x,y\right)  + P'_{N}\left(\frac{\partial^{I} f_{1}}{\partial z^{I}}\left(v_{1}\left(z'\right),x,y\right),\dots,\frac{\partial^{I} f_{N}}{\partial z^{I}}\left(v_{1}\left(z'\right),x,y\right)\right) =\chi_{N}'\left(v_{1}\left(z'\right),x,y\right),  \end{split} \label{222CCbu1}\end{equation}

 We apply Proposition $4.2$ from Mir\cite{Mir1} in (\ref{222CCbu1}).   Then, (\ref{MN1}) holds for $n\in\mathbb{N}$ respecting an induction procedure.

 \section{Analytic Systems in Codimension $2$}
 Let $M,M'\subset\mathbb{C}^{N+1}$ be two  Real-Analytic Submanifolds   as in  (\ref{W}) such that it exists a formal mapping  like (\ref{mapA}). We obtain
\begin{equation}   \displaystyle\sum_{ k\in\mathbb{N}}g_{k}(z)w^{k} = \left(Q_{2}+\varphi'\right)\left( \displaystyle\sum_{ k\in\mathbb{N}}f_{k}(z)w^{k} , \overline{ \displaystyle\sum_{ k\in\mathbb{N}}f_{k}(z)w^{k}}\right).   \label{V1bu11}\end{equation}
 \subsection{Case $Q_{1}\equiv\overline{Q}_{1}$ and  $Q_{2} \equiv\overline{Q}_{2}$} 
 Separating the real parts from the imaginary parts in (\ref{V1bu11}), we obtain
\begin{equation}  \frac{ \displaystyle\sum_{ k\in\mathbb{N}}\left(g_{k}(z)Z_{k}(x,z,\overline{z})-\overline{g_{k}(z)Z_{J}(x,z,\overline{z})}\right)}{2\sqrt{-1}} =  \left(\Im\varphi'\right)\left(\displaystyle\sum_{ k\in\mathbb{N}}f_{k}(z)Z_{k}(x,z,\overline{z}), \overline{ \displaystyle\sum_{ k\in\mathbb{N}}f_{k}(z)Z_{k}(x,z,\overline{z})} \right), \label{V21bu119} \end{equation}
 when (\ref{98q}) holds, and respectively, we obtain
\begin{equation}  \frac{ \displaystyle\sum_{ k\in\mathbb{N}}\left(g_{k}(z)Z_{k}(x,z,\overline{z})+\overline{g_{k}(z)Z_{J}(x,z,\overline{z})}\right)}{2 } =  \left( Q_{2}+\Re\varphi'\right)\left(\displaystyle\sum_{ k\in\mathbb{N}}f_{k}(z)Z_{k}(x,z,\overline{z}), \overline{ \displaystyle\sum_{ k\in\mathbb{N}}f_{k}(z)Z_{k}(x,z,\overline{z})} \right), \label{V22bu119} \end{equation}
 when (\ref{98q}) holds.   

The equations   (\ref{V21bu119}) and  (\ref{V22bu119}) are used   for suitable $\delta,\epsilon>0$ provided by Lemma \ref{tele21}. In particular, we chose by (\ref{98}) a point $z_{x}\in\mathbb{C}^{N}$ satisfying (\ref{61}).  Then, we    use  the Real-Analytic Submanifold $M_{x}$, which are resulted by (\ref{44}) and (\ref{441}) from The Implicit Function Theorem.

Next, we identify the coefficients of $x$ in   (\ref{V21bu119}) and  (\ref{V22bu119}).  We obtain that (\ref{MN2}) holds for $n=0$ and   the following analytic system:
\begin{equation} \begin{split}&\hspace{0.18 cm}  f_{1}\left(v_{1}\left(z'\right),x\right)+  P_{1}\left(f_{1}\left(v_{1}\left(z'\right),x\right),\dots,f_{N}\left(v_{1}\left(z'\right),x\right)\right)  =\chi_{1}\left(v_{1}\left(z'\right),x\right),\\&\hspace{0.18 cm} f_{2}\left(v_{1}\left(z'\right),x\right)   + P_{2}\left(f_{1}\left(v_{1}\left(z'\right),x\right),\dots,f_{N}\left(v_{1}\left(z'\right),x\right)\right) =\chi_{2}\left(v_{1}\left(z'\right),x\right), \\& \dots\dots\dots\dots\dots\dots\dots\dots\dots\dots\dots\dots\dots\dots\dots\dots\dots\dots\dots\dots\dots\dots\dots\dots\dots    \\&    f_{N}\left(v_{1}\left(z'\right),x\right)  +P_{N}\left(f_{1}\left(v_{1}\left(z'\right),x\right),\dots,f_{N}\left(v_{1}\left(z'\right),x\right)\right) =\chi_{N}\left(v_{1}\left(z'\right),x\right).  \end{split} \label{V6bu119}\end{equation}  
  
 We apply Proposition $4.2$ from Mir\cite{Mir1} in (\ref{V6}). Then, (\ref{MN1}) holds for $n=0$.

Next, we can take derivatives with respect to $z_{1},z_{2},\dots,z_{N}$ in (\ref{V22bu119}) of length $\left|I\right|=n$, for $I\in\mathbb{N}^{N}$  and a fixed $n\in\mathbb{N}$.  Repeating the previous procedure evaluating the resulted identity using (\ref{NX}), we conclude  the following analytic system
\begin{equation} \begin{split}&\hspace{0.17 cm} \frac{\partial^{I} f_{1}}{\partial z^{I}}\left(v_{1}\left(z'\right),x\right) + P'_{1}\left(\frac{\partial^{I} f_{1}}{\partial z^{I}}\left(v_{1}\left(z'\right),x\right),\dots,\frac{\partial^{I} f_{N}}{\partial z^{I}}\left(v_{1}\left(z'\right),x,y\right)\right) =\chi'_{1}\left(v_{1}\left(z'\right),x\right),  \\&\quad\hspace{0.15 cm} \frac{\partial^{I} f_{2}}{\partial z^{I}}\left(v_{1}\left(z'\right),x\right)  + P'_{2}\left(\frac{\partial^{I} f_{1}}{\partial z^{I}}\left(v_{1}(z),x\right),\dots,\frac{\partial^{I} f_{N}}{\partial z^{I}}\left(v_{1}\left(z'\right),x\right)\right)  =\chi'_{2}\left(v_{1}(z),x,y\right),   \\&\dots\dots\dots\dots\dots\dots\dots\dots\dots\dots\dots\dots\dots\dots\dots\dots\dots\dots\dots\dots\dots\dots\dots\dots\dots\dots\dots\dots\dots \\&\frac{\partial^{I}f_{N}}{\partial z^{I}}\left(v_{1}\left(z'\right),x,y\right)  + P'_{N}\left(\frac{\partial^{I} f_{1}}{\partial z^{I}}\left(v_{1}\left(z'\right),x,y\right),\dots,\frac{\partial^{I} f_{N}}{\partial z^{I}}\left(v_{1}\left(z'\right),x,y\right)\right) =\chi_{N}'\left(v_{1}\left(z'\right),x\right),  \end{split} \label{222CCbu119}\end{equation}

 We apply Proposition $4.2$ from Mir\cite{Mir1} in (\ref{222CCbu1}).   Then, (\ref{MN1}) holds for $n\in\mathbb{N}$ respecting an induction procedure.

 \subsection{Case $Q_{1}\not\equiv\overline{Q}_{1}$ and  $Q_{2}\not\equiv\overline{Q}_{2}$}
Separating the real parts from the imaginary parts in (\ref{V1bu11}), we obtain
\begin{equation}  \frac{ \displaystyle\sum_{ k\in\mathbb{N}}\left(g_{k}(z)Z_{k}(x,y,z,\overline{z})-\overline{g_{k}(z)Z_{J}(x,y,z,\overline{z})}\right)}{2\sqrt{-1}} =  \left(\Im Q_{2}+ \Im\varphi'\right)\left(\displaystyle\sum_{ k\in\mathbb{N}}f_{k}(z)Z_{k}(x,y,z,\overline{z}), \overline{ \displaystyle\sum_{ k\in\mathbb{N}}f_{k}(z)Z_{k}(x,y,z,\overline{z})} \right), \label{V21bu115} \end{equation}
 when (\ref{98q}) holds, and respectively, we obtain
\begin{equation}  \frac{ \displaystyle\sum_{ k\in\mathbb{N}}\left(g_{k}(z)Z_{k}(x,y,z,\overline{z})+\overline{g_{k}(z)Z_{J}(x,y,z,\overline{z})}\right)}{2 } =  \left(\Re Q_{2}+\Re\varphi'\right)\left(\displaystyle\sum_{ k\in\mathbb{N}}f_{k}(z)Z_{k}(x,y,z,\overline{z}), \overline{ \displaystyle\sum_{ k\in\mathbb{N}}f_{k}(z)Z_{k}(x,y,z,\overline{z})} \right), \label{V22bu115} \end{equation}
 when (\ref{98q}) holds.   

The equations   (\ref{V21bu115}) and  (\ref{V22bu115}) are used   for suitable $\delta,\epsilon>0$ provided by Lemma \ref{tele21}. In particular, we chose by (\ref{98q}) a point $z_{x}\in\mathbb{C}^{N}$ satisfying (\ref{61q}).  Then, we    use  the Real-Analytic Submanifold $M_{x}$, which are resulted by (\ref{44q}) and (\ref{441q}) from The Implicit Function Theorem.

Next, we identify the coefficients of $x$ in  (\ref{V21bu115}) and  (\ref{V22bu115}).  We obtain that (\ref{MN2}) holds for $n=0$ and   the following analytic system:
\begin{equation} \begin{split}&\hspace{0.18 cm}  f_{1}\left(v_{1}\left(z'\right),x,y\right)+  P_{1}\left(f_{1}\left(v_{1}\left(z'\right),x,y\right),\dots,f_{N}\left(v_{1}\left(z'\right),x,y\right)\right)  =\chi_{1}\left(v_{1}\left(z'\right),x,y\right),\\&\hspace{0.18 cm} f_{2}\left(v_{1}\left(z'\right),x,y\right)   + P_{2}\left(f_{1}\left(v_{1}\left(z'\right),x,y\right),\dots,f_{N}\left(v_{1}\left(z'\right),x,y\right)\right) =\chi_{2}\left(v_{1}\left(z'\right),x,y\right), \\& \dots\dots\dots\dots\dots\dots\dots\dots\dots\dots\dots\dots\dots\dots\dots\dots\dots\dots\dots\dots\dots\dots\dots\dots\dots    \\&    f_{N}\left(v_{1}\left(z'\right),x,y\right)  +P_{N}\left(f_{1}\left(v_{1}\left(z'\right),x,y\right),\dots,f_{N}\left(v_{1}\left(z'\right),x,y\right)\right) =\chi_{N}\left(v_{1}\left(z'\right),x,y\right).  \end{split} \label{V6bu115}\end{equation}  
  
 We apply Proposition $4.2$ from Mir\cite{Mir1} in (\ref{V6bu115}). Then, (\ref{MN1}) holds for $n=0$.

Next, we can take derivatives with respect to $z_{1},z_{2},\dots,z_{N}$ in (\ref{V22bu115}) of length $\left|I\right|=n$, for $I\in\mathbb{N}^{N}$  and a fixed $n\in\mathbb{N}$.  Repeating the previous procedure evaluating the resulted identity using (\ref{NX}), we conclude  the following analytic system
\begin{equation} \begin{split}&\hspace{0.17 cm} \frac{\partial^{I} f_{1}}{\partial z^{I}}\left(v_{1}\left(z'\right),x,y\right) + P'_{1}\left(\frac{\partial^{I} f_{1}}{\partial z^{I}}\left(v_{1}\left(z'\right),x,y\right),\dots,\frac{\partial^{I} f_{N}}{\partial z^{I}}\left(v_{1}\left(z'\right),x,y\right)\right) =\chi'_{1}\left(v_{1}\left(z'\right),x,y\right),  \\&\quad\hspace{0.15 cm} \frac{\partial^{I} f_{2}}{\partial z^{I}}\left(v_{1}\left(z'\right),x,y\right)  + P'_{2}\left(\frac{\partial^{I} f_{1}}{\partial z^{I}}\left(v_{1}(z),x,y\right),\dots,\frac{\partial^{I} f_{N}}{\partial z^{I}}\left(v_{1}\left(z'\right),x,y\right)\right)  =\chi'_{2}\left(v_{1}(z),x,y\right),   \\&\dots\dots\dots\dots\dots\dots\dots\dots\dots\dots\dots\dots\dots\dots\dots\dots\dots\dots\dots\dots\dots\dots\dots\dots\dots\dots\dots\dots\dots \\&\frac{\partial^{I}f_{N}}{\partial z^{I}}\left(v_{1}\left(z'\right),x,y\right)  + P'_{N}\left(\frac{\partial^{I} f_{1}}{\partial z^{I}}\left(v_{1}\left(z'\right),x,y\right),\dots,\frac{\partial^{I} f_{N}}{\partial z^{I}}\left(v_{1}\left(z'\right),x,y\right)\right) =\chi_{N}'\left(v_{1}\left(z'\right),x,y\right),  \end{split} \label{222CCbu11}\end{equation}

 We apply Proposition $4.2$ from Mir\cite{Mir1} in (\ref{222CCbu1}).   Then, (\ref{MN1}) holds for $n\in\mathbb{N}$ respecting an induction procedure.
  
 \subsection{Case $Q_{1}\not\equiv \overline{Q}_{1}$ and  $Q_{2}\equiv\overline{Q}_{2}$}
 
 Separating the real parts from the imaginary parts in (\ref{V1bu11}), we obtain
\begin{equation}  \frac{ \displaystyle\sum_{ k\in\mathbb{N}}\left(g_{k}(z)Z_{k}(x,y,z,\overline{z})-\overline{g_{k}(z)Z_{J}(x,y,z,\overline{z})}\right)}{2\sqrt{-1}} =  \left(  \Im\varphi'\right)\left(\displaystyle\sum_{ k\in\mathbb{N}}f_{k}(z)Z_{k}(x,y,z,\overline{z}), \overline{ \displaystyle\sum_{ k\in\mathbb{N}}f_{k}(z)Z_{k}(x,y,z,\overline{z})} \right), \label{V31} \end{equation}
 when (\ref{98q}) holds, and respectively, we obtain
\begin{equation}  \frac{ \displaystyle\sum_{ k\in\mathbb{N}}\left(g_{k}(z)Z_{k}(x,y,z,\overline{z})+\overline{g_{k}(z)Z_{J}(x,y,z,\overline{z})}\right)}{2 } =  \left(  Q_{2}+\Re\varphi'\right)\left(\displaystyle\sum_{ k\in\mathbb{N}}f_{k}(z)Z_{k}(x,y,z,\overline{z}), \overline{ \displaystyle\sum_{ k\in\mathbb{N}}f_{k}(z)Z_{k}(x,y,z,\overline{z})} \right), \label{V32} \end{equation}
 when (\ref{98q}) holds.   

The equations   (\ref{V31}) and  (\ref{V32}) are used   for suitable $\delta,\epsilon>0$ provided by Lemma \ref{tele21}. In particular, we chose by (\ref{98q}) a point $z_{x}\in\mathbb{C}^{N}$ satisfying (\ref{61q}).  Then, we    use  the Real-Analytic Submanifold $M_{x}$, which are resulted by (\ref{44q}) and (\ref{441q}) from The Implicit Function Theorem.

Next, we identify the coefficients of $x$ in  (\ref{V31}) and  (\ref{V32}).  We obtain that (\ref{MN2}) holds for $n=0$ and   the following analytic system:
\begin{equation} \begin{split}&\hspace{0.18 cm}  f_{1}\left(v_{1}\left(z'\right),x,y\right)+  P_{1}\left(f_{1}\left(v_{1}\left(z'\right),x,y\right),\dots,f_{N}\left(v_{1}\left(z'\right),x,y\right)\right)  =\chi_{1}\left(v_{1}\left(z'\right),x,y\right),\\&\hspace{0.18 cm} f_{2}\left(v_{1}\left(z'\right),x,y\right)   + P_{2}\left(f_{1}\left(v_{1}\left(z'\right),x,y\right),\dots,f_{N}\left(v_{1}\left(z'\right),x,y\right)\right) =\chi_{2}\left(v_{1}\left(z'\right),x,y\right), \\& \dots\dots\dots\dots\dots\dots\dots\dots\dots\dots\dots\dots\dots\dots\dots\dots\dots\dots\dots\dots\dots\dots\dots\dots\dots    \\&    f_{N}\left(v_{1}\left(z'\right),x,y\right)  +P_{N}\left(f_{1}\left(v_{1}\left(z'\right),x,y\right),\dots,f_{N}\left(v_{1}\left(z'\right),x,y\right)\right) =\chi_{N}\left(v_{1}\left(z'\right),x,y\right).  \end{split} \label{V65}\end{equation}  
  
 We apply Proposition $4.2$ from Mir\cite{Mir1} in (\ref{V65}). Then, (\ref{MN1}) holds for $n=0$.

Next, we can take derivatives with respect to $z_{1},z_{2},\dots,z_{N}$ in (\ref{V22bu115}) of length $\left|I\right|=n$, for $I\in\mathbb{N}^{N}$  and a fixed $n\in\mathbb{N}$.  Repeating the previous procedure evaluating the resulted identity using (\ref{NX}), we conclude  the following analytic system
\begin{equation} \begin{split}&\hspace{0.17 cm} \frac{\partial^{I} f_{1}}{\partial z^{I}}\left(v_{1}\left(z'\right),x,y\right) + P'_{1}\left(\frac{\partial^{I} f_{1}}{\partial z^{I}}\left(v_{1}\left(z'\right),x,y\right),\dots,\frac{\partial^{I} f_{N}}{\partial z^{I}}\left(v_{1}\left(z'\right),x,y\right)\right) =\chi'_{1}\left(v_{1}\left(z'\right),x,y\right),  \\&\quad\hspace{0.15 cm} \frac{\partial^{I} f_{2}}{\partial z^{I}}\left(v_{1}\left(z'\right),x,y\right)  + P'_{2}\left(\frac{\partial^{I} f_{1}}{\partial z^{I}}\left(v_{1}(z),x,y\right),\dots,\frac{\partial^{I} f_{N}}{\partial z^{I}}\left(v_{1}\left(z'\right),x,y\right)\right)  =\chi'_{2}\left(v_{1}(z),x,y\right),   \\&\dots\dots\dots\dots\dots\dots\dots\dots\dots\dots\dots\dots\dots\dots\dots\dots\dots\dots\dots\dots\dots\dots\dots\dots\dots\dots\dots\dots\dots \\&\frac{\partial^{I}f_{N}}{\partial z^{I}}\left(v_{1}\left(z'\right),x,y\right)  + P'_{N}\left(\frac{\partial^{I} f_{1}}{\partial z^{I}}\left(v_{1}\left(z'\right),x,y\right),\dots,\frac{\partial^{I} f_{N}}{\partial z^{I}}\left(v_{1}\left(z'\right),x,y\right)\right) =\chi_{N}'\left(v_{1}\left(z'\right),x,y\right),  \end{split} \label{29h}\end{equation}

 We apply Proposition $4.2$ from Mir\cite{Mir1} in (\ref{29h}).   Then, (\ref{MN1}) holds for $n\in\mathbb{N}$ respecting an induction procedure.
  
 \subsection{Case $Q_{1}\equiv\overline{Q}_{1}$ and  $Q_{2}\not\equiv\overline{Q}_{2}$} 

Separating the real parts from the imaginary parts in (\ref{V1bu11}), we obtain
\begin{equation}  \frac{ \displaystyle\sum_{ k\in\mathbb{N}}\left(g_{k}(z)Z_{k}(x,y,z,\overline{z})-\overline{g_{k}(z)Z_{J}(x,z,\overline{z})}\right)}{2\sqrt{-1}} =  \left(\Im Q_{2}+ \Im\varphi'\right)\left(\displaystyle\sum_{ k\in\mathbb{N}}f_{k}(z)Z_{k}(x,z,\overline{z}), \overline{ \displaystyle\sum_{ k\in\mathbb{N}}f_{k}(z)Z_{k}(x,z,\overline{z})} \right), \label{Ve1} \end{equation}
 when (\ref{98}) holds, and respectively, we obtain
\begin{equation}  \frac{ \displaystyle\sum_{ k\in\mathbb{N}}\left(g_{k}(z)Z_{k}(x,y,z,\overline{z})+\overline{g_{k}(z)Z_{J}(x,y,z,\overline{z})}\right)}{2 } =  \left(\Re Q_{2}+\Re\varphi'\right)\left(\displaystyle\sum_{ k\in\mathbb{N}}f_{k}(z)Z_{k}(x,y,z,\overline{z}), \overline{ \displaystyle\sum_{ k\in\mathbb{N}}f_{k}(z)Z_{k}(x,y,z,\overline{z})} \right), \label{Ve2} \end{equation}
 when (\ref{98q}) holds.   

The equations   (\ref{Ve1}) and  (\ref{Ve2}) are used   for suitable $\delta,\epsilon>0$ provided by Lemma \ref{tele21}. In particular, we chose by (\ref{98}) a point $z_{x}\in\mathbb{C}^{N}$ satisfying (\ref{61}).  Then, we    use  the Real-Analytic Submanifold $M_{x}$, which are resulted by (\ref{44}) and (\ref{441}) from The Implicit Function Theorem.

Next, we identify the coefficients of $x$ in  (\ref{Ve1}) and  (\ref{Ve2}).  We obtain that (\ref{MN2}) holds for $n=0$ and   the following analytic system:
\begin{equation} \begin{split}&\hspace{0.18 cm}  f_{1}\left(v_{1}\left(z'\right),x\right)+  P_{1}\left(f_{1}\left(v_{1}\left(z'\right),x\right),\dots,f_{N}\left(v_{1}\left(z'\right),x,y\right)\right)  =\chi_{1}\left(v_{1}\left(z'\right),x\right),\\&\hspace{0.18 cm} f_{2}\left(v_{1}\left(z'\right),x\right)   + P_{2}\left(f_{1}\left(v_{1}\left(z'\right),x\right),\dots,f_{N}\left(v_{1}\left(z'\right),x\right)\right) =\chi_{2}\left(v_{1}\left(z'\right),x\right), \\& \dots\dots\dots\dots\dots\dots\dots\dots\dots\dots\dots\dots\dots\dots\dots\dots\dots\dots\dots\dots\dots\dots\dots\dots\dots    \\&    f_{N}\left(v_{1}\left(z'\right),x\right)  +P_{N}\left(f_{1}\left(v_{1}\left(z'\right),x\right),\dots,f_{N}\left(v_{1}\left(z'\right),x\right)\right) =\chi_{N}\left(v_{1}\left(z'\right),x\right).  \end{split} \label{bu}\end{equation}  
  
 We apply Proposition $4.2$ from Mir\cite{Mir1} in (\ref{bu}). Then, (\ref{MN1}) holds for $n=0$.

Next, we can take derivatives with respect to $z_{1},z_{2},\dots,z_{N}$ in (\ref{V22bu115}) of length $\left|I\right|=n$, for $I\in\mathbb{N}^{N}$  and a fixed $n\in\mathbb{N}$.  Repeating the previous procedure evaluating the resulted identity using (\ref{NX}), we conclude  the following analytic system
\begin{equation} \begin{split}&\hspace{0.17 cm} \frac{\partial^{I} f_{1}}{\partial z^{I}}\left(v_{1}\left(z'\right),x\right) + P'_{1}\left(\frac{\partial^{I} f_{1}}{\partial z^{I}}\left(v_{1}\left(z'\right),x\right),\dots,\frac{\partial^{I} f_{N}}{\partial z^{I}}\left(v_{1}\left(z'\right),x\right)\right) =\chi'_{1}\left(v_{1}\left(z'\right),x\right),  \\&\quad\hspace{0.15 cm} \frac{\partial^{I} f_{2}}{\partial z^{I}}\left(v_{1}\left(z'\right),x\right)  + P'_{2}\left(\frac{\partial^{I} f_{1}}{\partial z^{I}}\left(v_{1}(z),x\right),\dots,\frac{\partial^{I} f_{N}}{\partial z^{I}}\left(v_{1}\left(z'\right),x\right)\right)  =\chi'_{2}\left(v_{1}(z),x,y\right),   \\&\dots\dots\dots\dots\dots\dots\dots\dots\dots\dots\dots\dots\dots\dots\dots\dots\dots\dots\dots\dots\dots\dots\dots\dots\dots\dots\dots\dots\dots \\&\frac{\partial^{I}f_{N}}{\partial z^{I}}\left(v_{1}\left(z'\right),x\right)  + P'_{N}\left(\frac{\partial^{I} f_{1}}{\partial z^{I}}\left(v_{1}\left(z'\right),x\right),\dots,\frac{\partial^{I} f_{N}}{\partial z^{I}}\left(v_{1}\left(z'\right),x\right)\right) =\chi_{N}'\left(v_{1}\left(z'\right),x\right),  \end{split} \label{5a}\end{equation}

 We apply Proposition $4.2$ from Mir\cite{Mir1} in (\ref{5a}).   Then, (\ref{MN1}) holds for $n\in\mathbb{N}$ respecting an induction procedure.

  \section{Analytic Systems: CR Singularities and Real-Hypersurfaces in $\mathbb{C}^{N'+1}$}
 Let $M\subset\mathbb{C}^{N+1}$ be a  Real-Analytic Submanifold    as in  (\ref{W}) and $M'\subset\mathbb{C}^{N'+1}$ a Real-Analytic Hypersurface such that it exists a formal mapping  like (\ref{mapA}). We obtain
\begin{equation}\begin{split}& \frac{  \displaystyle\sum_{ k\in\mathbb{N}}g_{k}(z)w^{k}-\overline{\displaystyle\sum_{ k\in\mathbb{N}}g_{k}(z)w^{k}}}{2\sqrt{-1}} = \left(\frac{  \displaystyle\sum_{ k\in\mathbb{N}}g_{k}(z)w^{k}+\overline{\displaystyle\sum_{ k\in\mathbb{N}}g_{k}(z)w^{k}}}{2}\right)^{m}P\left(\displaystyle\sum_{ k\in\mathbb{N}}f_{k}(z)w^{k}  ,\right.\\&\left.\quad\quad \overline{ \displaystyle\sum_{ k\in\mathbb{N}}f_{k}(z)w^{k}}\right)+\varphi'\left( \displaystyle\sum_{ k\in\mathbb{N}}f_{k}(z)w^{k} , \overline{ \displaystyle\sum_{ k\in\mathbb{N}}f_{k}(z)w^{k}},\frac{  \displaystyle\sum_{ k\in\mathbb{N}}g_{k}(z)w^{k}+\overline{\displaystyle\sum_{ k\in\mathbb{N}}g_{k}(z)w^{k}}}{2}\right). \end{split}  \label{W}\end{equation}
\subsection{Case $Q\not\equiv \overline{Q}$}  Separating the real parts from the imaginary parts in (\ref{Ww}), we obtain
\begin{equation}\begin{split}& \quad\quad\frac{  \displaystyle\sum_{ k\in\mathbb{N}}g_{k}(z)Z_{k}(x,y,z,\overline{z})-\overline{\displaystyle\sum_{ k\in\mathbb{N}}g_{k}(z)Z_{k}(x,y,z,\overline{z})}}{2\sqrt{-1}} = \left(\frac{  \displaystyle\sum_{ k\in\mathbb{N}}g_{k}(z)w^{k}+\overline{\displaystyle\sum_{ k\in\mathbb{N}}g_{k}(z)w^{k}}}{2}\right)^{m}P\left(\displaystyle\sum_{ k\in\mathbb{N}}f_{k}(z)Z_{k}(x,y,z,\overline{z})  ,\right.\\&\left. \overline{ \displaystyle\sum_{ k\in\mathbb{N}}f_{k}(z)Z_{k}(x,y,z,\overline{z})}\right)+\varphi'\left( \displaystyle\sum_{ k\in\mathbb{N}}f_{k}(z)Z_{k}(x,y,z,\overline{z}) , \overline{ \displaystyle\sum_{ k\in\mathbb{N}}f_{k}(z)Z_{k}(x,y,z,\overline{z})},\frac{  \displaystyle\sum_{ k\in\mathbb{N}}g_{k}(z)w^{k}+\overline{\displaystyle\sum_{ k\in\mathbb{N}}g_{k}(z)Z_{k}(x,y,z,\overline{z})}}{2}\right), \end{split}  \label{V311a}\end{equation}
when (\ref{98q}) holds.

The equation   (\ref{V311a}) is used   for suitable $\delta,\epsilon>0$ provided by Lemma \ref{tele21}. In particular, we chose by (\ref{98q}) a point $z_{x}\in\mathbb{C}^{N}$ satisfying (\ref{61q}).  Then, we    use  the Real-Analytic Submanifold $M_{x}$, which are resulted by (\ref{44q}) and (\ref{441q}) from The Implicit Function Theorem.

Next, we identify the coefficients of $x$ in  (\ref{V311a}).  We obtain that (\ref{MN2}) holds for $n=0$ and   the following analytic system:
\begin{equation} \begin{split}&\hspace{0.18 cm}  f_{1}\left(v_{1}\left(z'\right),x,y\right)+  P_{1}\left(f_{1}\left(v_{1}\left(z'\right),x,y\right),\dots,f_{N}\left(v_{1}\left(z'\right),x,y\right)\right)  =\chi_{1}\left(v_{1}\left(z'\right),x,y\right),\\&\hspace{0.18 cm} f_{2}\left(v_{1}\left(z'\right),x,y\right)   + P_{2}\left(f_{1}\left(v_{1}\left(z'\right),x,y\right),\dots,f_{N}\left(v_{1}\left(z'\right),x,y\right)\right) =\chi_{2}\left(v_{1}\left(z'\right),x,y\right), \\& \dots\dots\dots\dots\dots\dots\dots\dots\dots\dots\dots\dots\dots\dots\dots\dots\dots\dots\dots\dots\dots\dots\dots\dots\dots    \\&    f_{N}\left(v_{1}\left(z'\right),x,y\right)  +P_{N}\left(f_{1}\left(v_{1}\left(z'\right),x,y\right),\dots,f_{N}\left(v_{1}\left(z'\right),x,y\right)\right) =\chi_{N}\left(v_{1}\left(z'\right),x,y\right).  \end{split} \label{V651y}\end{equation}  
  
 We apply Proposition $4.2$ from Mir\cite{Mir1} in (\ref{V651y}). Then, (\ref{MN1}) holds for $n=0$.

Next, we can take derivatives with respect to $z_{1},z_{2},\dots,z_{N}$ in (\ref{W}) of length $\left|I\right|=n$, for $I\in\mathbb{N}^{N}$  and a fixed $n\in\mathbb{N}$.  Repeating the previous procedure evaluating the resulted identity using (\ref{NX}), we conclude  the following analytic system
\begin{equation} \begin{split}&\hspace{0.17 cm} \frac{\partial^{I} f_{1}}{\partial z^{I}}\left(v_{1}\left(z'\right),x,y\right) + P'_{1}\left(\frac{\partial^{I} f_{1}}{\partial z^{I}}\left(v_{1}\left(z'\right),x,y\right),\dots,\frac{\partial^{I} f_{N}}{\partial z^{I}}\left(v_{1}\left(z'\right),x,y\right)\right) =\chi'_{1}\left(v_{1}\left(z'\right),x,y\right),  \\&\quad\hspace{0.15 cm} \frac{\partial^{I} f_{2}}{\partial z^{I}}\left(v_{1}\left(z'\right),x,y\right)  + P'_{2}\left(\frac{\partial^{I} f_{1}}{\partial z^{I}}\left(v_{1}(z),x,y\right),\dots,\frac{\partial^{I} f_{N}}{\partial z^{I}}\left(v_{1}\left(z'\right),x,y\right)\right)  =\chi'_{2}\left(v_{1}(z),x,y\right),   \\&\dots\dots\dots\dots\dots\dots\dots\dots\dots\dots\dots\dots\dots\dots\dots\dots\dots\dots\dots\dots\dots\dots\dots\dots\dots\dots\dots\dots\dots \\&\frac{\partial^{I}f_{N}}{\partial z^{I}}\left(v_{1}\left(z'\right),x,y\right)  + P'_{N}\left(\frac{\partial^{I} f_{1}}{\partial z^{I}}\left(v_{1}\left(z'\right),x,y\right),\dots,\frac{\partial^{I} f_{N}}{\partial z^{I}}\left(v_{1}\left(z'\right),x,y\right)\right) =\chi_{N}'\left(v_{1}\left(z'\right),x,y\right),  \end{split} \label{29h1}\end{equation}

 We apply Proposition $4.2$ from Mir\cite{Mir1} in (\ref{29h1}).   Then, (\ref{MN1}) holds for $n\in\mathbb{N}$ respecting an induction procedure.
  
 \subsection{Case $Q\equiv\overline{Q}$} Separating the real parts from the imaginary parts in (\ref{W}), we obtain
 \begin{equation}\begin{split}& \frac{  \displaystyle\sum_{ k\in\mathbb{N}}g_{k}(z)Z_{k}(x,z,\overline{z})-\overline{\displaystyle\sum_{ k\in\mathbb{N}}g_{k}(z)Z_{k}(x,z,\overline{z})}}{2\sqrt{-1}} = \left(\frac{  \displaystyle\sum_{ k\in\mathbb{N}}g_{k}(z)Z_{k}(x,z,\overline{z})+\overline{\displaystyle\sum_{ k\in\mathbb{N}}g_{k}(z)Z_{k}(x,z,\overline{z})}}{2}\right)^{m}P\left(\displaystyle\sum_{ k\in\mathbb{N}}f_{k}(z)Z_{k}(x,z,\overline{z})  ,\right.\\&\left.\quad\quad \overline{ \displaystyle\sum_{ k\in\mathbb{N}}f_{k}(z)Z_{k}(x,z,\overline{z})}\right)+\varphi'\left( \displaystyle\sum_{ k\in\mathbb{N}}f_{k}(z)Z_{k}(x,z,\overline{z}) , \overline{ \displaystyle\sum_{ k\in\mathbb{N}}f_{k}(z)Z_{k}(x,z,\overline{z})},\frac{  \displaystyle\sum_{ k\in\mathbb{N}}g_{k}(z)Z_{k}(x,z,\overline{z})+\overline{\displaystyle\sum_{ k\in\mathbb{N}}g_{k}(z)Z_{k}(x,z,\overline{z})}}{2}\right). \end{split}  \label{Ve1ab}\end{equation}
when (\ref{98}) holds. 

The equation   (\ref{Ve1ab}) is used   for suitable $\delta,\epsilon>0$ provided by Lemma \ref{tele21}. In particular, we chose by (\ref{98}) a point $z_{x}\in\mathbb{C}^{N}$ satisfying (\ref{61}).  Then, we    use  the Real-Analytic Submanifold $M_{x}$, which are resulted by (\ref{44}) and (\ref{441}) from The Implicit Function Theorem.

Next, we identify the coefficients of $x$ in  (\ref{Ve1ab}).  We obtain that (\ref{MN2}) holds for $n=0$ and   the following analytic system:
\begin{equation} \begin{split}&\hspace{0.18 cm}  f_{1}\left(v_{1}\left(z'\right),x\right)+  P_{1}\left(f_{1}\left(v_{1}\left(z'\right),x\right),\dots,f_{N}\left(v_{1}\left(z'\right),x,y\right)\right)  =\chi_{1}\left(v_{1}\left(z'\right),x\right),\\&\hspace{0.18 cm} f_{2}\left(v_{1}\left(z'\right),x\right)   + P_{2}\left(f_{1}\left(v_{1}\left(z'\right),x\right),\dots,f_{N}\left(v_{1}\left(z'\right),x\right)\right) =\chi_{2}\left(v_{1}\left(z'\right),x\right), \\& \dots\dots\dots\dots\dots\dots\dots\dots\dots\dots\dots\dots\dots\dots\dots\dots\dots\dots\dots\dots\dots\dots\dots\dots\dots    \\&    f_{N}\left(v_{1}\left(z'\right),x\right)  +P_{N}\left(f_{1}\left(v_{1}\left(z'\right),x\right),\dots,f_{N}\left(v_{1}\left(z'\right),x\right)\right) =\chi_{N}\left(v_{1}\left(z'\right),x\right).  \end{split} \label{buaab}\end{equation}  
  
 We apply Proposition $4.2$ from Mir\cite{Mir1} in (\ref{buaab}). Then, (\ref{MN1}) holds for $n=0$.

Next, we can take derivatives with respect to $z_{1},z_{2},\dots,z_{N}$ in (\ref{W}) of length $\left|I\right|=n$, for $I\in\mathbb{N}^{N}$  and a fixed $n\in\mathbb{N}$.  Repeating the previous procedure evaluating the resulted identity using (\ref{NX}), we conclude  the following analytic system
\begin{equation} \begin{split}&\hspace{0.17 cm} \frac{\partial^{I} f_{1}}{\partial z^{I}}\left(v_{1}\left(z'\right),x\right) + P'_{1}\left(\frac{\partial^{I} f_{1}}{\partial z^{I}}\left(v_{1}\left(z'\right),x\right),\dots,\frac{\partial^{I} f_{N}}{\partial z^{I}}\left(v_{1}\left(z'\right),x\right)\right) =\chi'_{1}\left(v_{1}\left(z'\right),x\right),  \\&\quad\hspace{0.15 cm} \frac{\partial^{I} f_{2}}{\partial z^{I}}\left(v_{1}\left(z'\right),x\right)  + P'_{2}\left(\frac{\partial^{I} f_{1}}{\partial z^{I}}\left(v_{1}(z),x\right),\dots,\frac{\partial^{I} f_{N}}{\partial z^{I}}\left(v_{1}\left(z'\right),x\right)\right)  =\chi'_{2}\left(v_{1}(z),x,y\right),   \\&\dots\dots\dots\dots\dots\dots\dots\dots\dots\dots\dots\dots\dots\dots\dots\dots\dots\dots\dots\dots\dots\dots\dots\dots\dots\dots\dots\dots\dots \\&\frac{\partial^{I}f_{N}}{\partial z^{I}}\left(v_{1}\left(z'\right),x\right)  + P'_{N}\left(\frac{\partial^{I} f_{1}}{\partial z^{I}}\left(v_{1}\left(z'\right),x\right),\dots,\frac{\partial^{I} f_{N}}{\partial z^{I}}\left(v_{1}\left(z'\right),x\right)\right) =\chi_{N}'\left(v_{1}\left(z'\right),x\right),  \end{split} \label{5aaab}\end{equation}

 We apply Proposition $4.2$ from Mir\cite{Mir1} in (\ref{5aaab}).   Then, (\ref{MN1}) holds for $n\in\mathbb{N}$ respecting an induction procedure.

  \section{Analytic Systems: CR Singularities and Real-Surfaces in $\mathbb{C}^{2}$}
 Let $M\subset\mathbb{C}^{N+1}$ be a  Real-Analytic Submanifold    as in  (\ref{W}) and $M'\subset\mathbb{C}^{2}$ a Real-Analytic Surface such that it exists a formal mapping  like (\ref{mapA}). We obtain
\begin{equation}   \displaystyle\sum_{ k\in\mathbb{N}}g_{k}(z)w^{k} = \left(P+\varphi'\right)\left( \displaystyle\sum_{ k\in\mathbb{N}}f_{k}(z)w^{k} , \overline{ \displaystyle\sum_{ k\in\mathbb{N}}f_{k}(z)w^{k}}\right).   \label{Ww}\end{equation}
 \subsection{Case $Q\not\equiv \overline{Q}$}  Separating the real parts from the imaginary parts in (\ref{Ww}), we obtain
\begin{equation}  \frac{ \displaystyle\sum_{ k\in\mathbb{N}}\left(g_{k}(z)Z_{k}(x,y,z,\overline{z})-\overline{g_{k}(z)Z_{J}(x,y,z,\overline{z})}\right)}{2\sqrt{-1}} =  \left( \Im P+ \Im\varphi'\right)\left(\displaystyle\sum_{ k\in\mathbb{N}}f_{k}(z)Z_{k}(x,y,z,\overline{z}), \overline{ \displaystyle\sum_{ k\in\mathbb{N}}f_{k}(z)Z_{k}(x,y,z,\overline{z})} \right), \label{V311} \end{equation}
 when (\ref{98q}) holds, and respectively, we obtain
\begin{equation}  \frac{ \displaystyle\sum_{ k\in\mathbb{N}}\left(g_{k}(z)Z_{k}(x,y,z,\overline{z})+\overline{g_{k}(z)Z_{J}(x,y,z,\overline{z})}\right)}{2 } =  \left( \Re P+\Re\varphi'\right)\left(\displaystyle\sum_{ k\in\mathbb{N}}f_{k}(z)Z_{k}(x,y,z,\overline{z}), \overline{ \displaystyle\sum_{ k\in\mathbb{N}}f_{k}(z)Z_{k}(x,y,z,\overline{z})} \right), \label{V321} \end{equation}
 when (\ref{98q}) holds.   

The equations   (\ref{V311}) and  (\ref{V321}) are used   for suitable $\delta,\epsilon>0$ provided by Lemma \ref{tele21}. In particular, we chose by (\ref{98q}) a point $z_{x}\in\mathbb{C}^{N}$ satisfying (\ref{61q}).  Then, we    use  the Real-Analytic Submanifold $M_{x}$, which are resulted by (\ref{44q}) and (\ref{441q}) from The Implicit Function Theorem.

Next, we identify the coefficients of $x$ in  (\ref{V311}) and  (\ref{V321}).  We obtain that (\ref{MN2}) holds for $n=0$ and   the following analytic system:
\begin{equation} \begin{split}&\hspace{0.18 cm}  f_{1}\left(v_{1}\left(z'\right),x,y\right)+  P_{1}\left(f_{1}\left(v_{1}\left(z'\right),x,y\right),\dots,f_{N}\left(v_{1}\left(z'\right),x,y\right)\right)  =\chi_{1}\left(v_{1}\left(z'\right),x,y\right),\\&\hspace{0.18 cm} f_{2}\left(v_{1}\left(z'\right),x,y\right)   + P_{2}\left(f_{1}\left(v_{1}\left(z'\right),x,y\right),\dots,f_{N}\left(v_{1}\left(z'\right),x,y\right)\right) =\chi_{2}\left(v_{1}\left(z'\right),x,y\right), \\& \dots\dots\dots\dots\dots\dots\dots\dots\dots\dots\dots\dots\dots\dots\dots\dots\dots\dots\dots\dots\dots\dots\dots\dots\dots    \\&    f_{N}\left(v_{1}\left(z'\right),x,y\right)  +P_{N}\left(f_{1}\left(v_{1}\left(z'\right),x,y\right),\dots,f_{N}\left(v_{1}\left(z'\right),x,y\right)\right) =\chi_{N}\left(v_{1}\left(z'\right),x,y\right).  \end{split} \label{V651}\end{equation}  
  
 We apply Proposition $4.2$ from Mir\cite{Mir1} in (\ref{V651}). Then, (\ref{MN1}) holds for $n=0$.

Next, we can take derivatives with respect to $z_{1},z_{2},\dots,z_{N}$ in (\ref{V22bu115}) of length $\left|I\right|=n$, for $I\in\mathbb{N}^{N}$  and a fixed $n\in\mathbb{N}$.  Repeating the previous procedure evaluating the resulted identity using (\ref{NX}), we conclude  the following analytic system
\begin{equation} \begin{split}&\hspace{0.17 cm} \frac{\partial^{I} f_{1}}{\partial z^{I}}\left(v_{1}\left(z'\right),x,y\right) + P'_{1}\left(\frac{\partial^{I} f_{1}}{\partial z^{I}}\left(v_{1}\left(z'\right),x,y\right),\dots,\frac{\partial^{I} f_{N}}{\partial z^{I}}\left(v_{1}\left(z'\right),x,y\right)\right) =\chi'_{1}\left(v_{1}\left(z'\right),x,y\right),  \\&\quad\hspace{0.15 cm} \frac{\partial^{I} f_{2}}{\partial z^{I}}\left(v_{1}\left(z'\right),x,y\right)  + P'_{2}\left(\frac{\partial^{I} f_{1}}{\partial z^{I}}\left(v_{1}(z),x,y\right),\dots,\frac{\partial^{I} f_{N}}{\partial z^{I}}\left(v_{1}\left(z'\right),x,y\right)\right)  =\chi'_{2}\left(v_{1}(z),x,y\right),   \\&\dots\dots\dots\dots\dots\dots\dots\dots\dots\dots\dots\dots\dots\dots\dots\dots\dots\dots\dots\dots\dots\dots\dots\dots\dots\dots\dots\dots\dots \\&\frac{\partial^{I}f_{N}}{\partial z^{I}}\left(v_{1}\left(z'\right),x,y\right)  + P'_{N}\left(\frac{\partial^{I} f_{1}}{\partial z^{I}}\left(v_{1}\left(z'\right),x,y\right),\dots,\frac{\partial^{I} f_{N}}{\partial z^{I}}\left(v_{1}\left(z'\right),x,y\right)\right) =\chi_{N}'\left(v_{1}\left(z'\right),x,y\right),  \end{split} \label{29h1}\end{equation}

 We apply Proposition $4.2$ from Mir\cite{Mir1} in (\ref{29h1}).   Then, (\ref{MN1}) holds for $n\in\mathbb{N}$ respecting an induction procedure.
  
 \subsection{Case $Q\equiv\overline{Q}$} Separating the real parts from the imaginary parts in (\ref{Ww}), we obtain
\begin{equation}  \frac{ \displaystyle\sum_{ k\in\mathbb{N}}\left(g_{k}(z)Z_{k}(x,z,\overline{z})-\overline{g_{k}(z)Z_{J}(x,z,\overline{z})}\right)}{2\sqrt{-1}} =  \left(\Im P+ \Im\varphi'\right)\left(\displaystyle\sum_{ k\in\mathbb{N}}f_{k}(z)Z_{k}(x,z,\overline{z}), \overline{ \displaystyle\sum_{ k\in\mathbb{N}}f_{k}(z)Z_{k}(x,z,\overline{z})} \right), \label{Ve1a} \end{equation}
 when (\ref{98}) holds, and respectively, we obtain
\begin{equation}  \frac{ \displaystyle\sum_{ k\in\mathbb{N}}\left(g_{k}(z)Z_{k}(x,z,\overline{z})+\overline{g_{k}(z)Z_{J}(x,z,\overline{z})}\right)}{2 } =  \left(\Re P+\Re\varphi'\right)\left(\displaystyle\sum_{ k\in\mathbb{N}}f_{k}(z)Z_{k}(x,z,\overline{z}), \overline{ \displaystyle\sum_{ k\in\mathbb{N}}f_{k}(z)Z_{k}(x,z,\overline{z})} \right), \label{Ve2a} \end{equation}
 when (\ref{98q}) holds.   

The equations   (\ref{Ve1a}) and  (\ref{Ve2a}) are used   for suitable $\delta,\epsilon>0$ provided by Lemma \ref{tele21}. In particular, we chose by (\ref{98}) a point $z_{x}\in\mathbb{C}^{N}$ satisfying (\ref{61}).  Then, we    use  the Real-Analytic Submanifold $M_{x}$, which are resulted by (\ref{44}) and (\ref{441}) from The Implicit Function Theorem.

Next, we identify the coefficients of $x$ in  (\ref{Ve1a}) and  (\ref{Ve2a}).  We obtain that (\ref{MN2}) holds for $n=0$ and   the following analytic system:
\begin{equation} \begin{split}&\hspace{0.18 cm}  f_{1}\left(v_{1}\left(z'\right),x\right)+  P_{1}\left(f_{1}\left(v_{1}\left(z'\right),x\right),\dots,f_{N}\left(v_{1}\left(z'\right),x,y\right)\right)  =\chi_{1}\left(v_{1}\left(z'\right),x\right),\\&\hspace{0.18 cm} f_{2}\left(v_{1}\left(z'\right),x\right)   + P_{2}\left(f_{1}\left(v_{1}\left(z'\right),x\right),\dots,f_{N}\left(v_{1}\left(z'\right),x\right)\right) =\chi_{2}\left(v_{1}\left(z'\right),x\right), \\& \dots\dots\dots\dots\dots\dots\dots\dots\dots\dots\dots\dots\dots\dots\dots\dots\dots\dots\dots\dots\dots\dots\dots\dots\dots    \\&    f_{N}\left(v_{1}\left(z'\right),x\right)  +P_{N}\left(f_{1}\left(v_{1}\left(z'\right),x\right),\dots,f_{N}\left(v_{1}\left(z'\right),x\right)\right) =\chi_{N}\left(v_{1}\left(z'\right),x\right).  \end{split} \label{buaa}\end{equation}  
  
 We apply Proposition $4.2$ from Mir\cite{Mir1} in (\ref{buaa}). Then, (\ref{MN1}) holds for $n=0$.

Next, we can take derivatives with respect to $z_{1},z_{2},\dots,z_{N}$ in (\ref{V22bu115}) of length $\left|I\right|=n$, for $I\in\mathbb{N}^{N}$  and a fixed $n\in\mathbb{N}$.  Repeating the previous procedure evaluating the resulted identity using (\ref{NX}), we conclude  the following analytic system
\begin{equation} \begin{split}&\hspace{0.17 cm} \frac{\partial^{I} f_{1}}{\partial z^{I}}\left(v_{1}\left(z'\right),x\right) + P'_{1}\left(\frac{\partial^{I} f_{1}}{\partial z^{I}}\left(v_{1}\left(z'\right),x\right),\dots,\frac{\partial^{I} f_{N}}{\partial z^{I}}\left(v_{1}\left(z'\right),x\right)\right) =\chi'_{1}\left(v_{1}\left(z'\right),x\right),  \\&\quad\hspace{0.15 cm} \frac{\partial^{I} f_{2}}{\partial z^{I}}\left(v_{1}\left(z'\right),x\right)  + P'_{2}\left(\frac{\partial^{I} f_{1}}{\partial z^{I}}\left(v_{1}(z),x\right),\dots,\frac{\partial^{I} f_{N}}{\partial z^{I}}\left(v_{1}\left(z'\right),x\right)\right)  =\chi'_{2}\left(v_{1}(z),x,y\right),   \\&\dots\dots\dots\dots\dots\dots\dots\dots\dots\dots\dots\dots\dots\dots\dots\dots\dots\dots\dots\dots\dots\dots\dots\dots\dots\dots\dots\dots\dots \\&\frac{\partial^{I}f_{N}}{\partial z^{I}}\left(v_{1}\left(z'\right),x\right)  + P'_{N}\left(\frac{\partial^{I} f_{1}}{\partial z^{I}}\left(v_{1}\left(z'\right),x\right),\dots,\frac{\partial^{I} f_{N}}{\partial z^{I}}\left(v_{1}\left(z'\right),x\right)\right) =\chi_{N}'\left(v_{1}\left(z'\right),x\right),  \end{split} \label{5aaa}\end{equation}

 We apply Proposition $4.2$ from Mir\cite{Mir1} in (\ref{5aaa}).   Then, (\ref{MN1}) holds for $n\in\mathbb{N}$ respecting an induction procedure.

\section{Proof of Theorem \ref{t1}}
We  look  at  $z=\left(z_{1},z_{2},\dots,z_{N}\right)$ and $\overline{z}=\left(\overline{z}_{1},\overline{z}_{2},\dots,\overline{z}_{N}\right)$ as independent variables, replacing  $\overline{z}$ with 
$  \zeta=\left(\zeta_{1},\zeta_{2},\dots,\zeta_{N}\right) $. Using $z'=\left(z_{1},z_{2},\dots,z_{N-1}\right)$ and $\zeta'=\left(\zeta_{1},\zeta_{2},\dots,\zeta_{N-1}\right),$ we consider the ring of the formal power series in $ x,t, z',z,\zeta$,  denoted by $\mathbb{C}\left[\left[   x,t,z',\zeta\right]\right]$, where $t$ is just the complexification of $x$. We assume $x,t$ near $0\in\mathbb{C}^{N}$ by an eventual translation.  We complexify the equations
(\ref{V21}) and (\ref{V22}). We apply the Approximation Theorem of Artin\cite{A}. We obtain the desired convergence after few non-trivial computations.

\section{Proof of Theorem \ref{t2}}
We  look  at  $z=\left(z_{1},z_{2},\dots,z_{N}\right)$ and $\overline{z}=\left(\overline{z}_{1},\overline{z}_{2},\dots,\overline{z}_{N}\right)$ as independent variables, replacing  $\overline{z}$ with 
$  \zeta=\left(\zeta_{1},\zeta_{2},\dots,\zeta_{N}\right) $. Using $z'=\left(z_{1},z_{2},\dots,z_{N-1}\right)$ and $\zeta'=\left(\zeta_{1},\zeta_{2},\dots,\zeta_{N-1}\right),$ we consider the ring of the formal power series in $ x,t, z',z,\zeta$,  denoted by $\mathbb{C}\left[\left[   x,t,z',\zeta\right]\right]$, where $t$ is just the complexification of $x$. We assume $x,t$ near $0\in\mathbb{C}^{N}$ by an eventual translation.  We complexify the equations
(\ref{V21bu1}) and (\ref{V22bu1}). We apply the Approximation Theorem of Artin\cite{A}. We obtain the desired convergence after few non-trivial computations.
  \section{Proof of Theorem \ref{t3}}
We  look  at  $z=\left(z_{1},z_{2},\dots,z_{N}\right)$ and $\overline{z}=\left(\overline{z}_{1},\overline{z}_{2},\dots,\overline{z}_{N}\right)$ as independent variables, replacing  $\overline{z}$ with 
$  \zeta=\left(\zeta_{1},\zeta_{2},\dots,\zeta_{N}\right) $. Using $z'=\left(z_{1},z_{2},\dots,z_{N-1}\right)$ and $\zeta'=\left(\zeta_{1},\zeta_{2},\dots,\zeta_{N-1}\right),$ we consider the ring of the formal power series in $ x,y,t,r, z',z,\zeta$,  denoted by $\mathbb{C}\left[\left[   x,t,z',\zeta\right]\right]$, where $t,r$ is just the complexification of $x,y$. We assume $x,y,t,r$ near $0\in\mathbb{C}^{N}$ by an eventual translation.  We complexify the equations
(\ref{V21bu119}) and (\ref{V22bu119}). We apply the Approximation Theorem of Artin\cite{A}. We obtain the desired convergence after few non-trivial computations.
\section{Proof of Theorem \ref{t4}}
We  look  at  $z=\left(z_{1},z_{2},\dots,z_{N}\right)$ and $\overline{z}=\left(\overline{z}_{1},\overline{z}_{2},\dots,\overline{z}_{N}\right)$ as independent variables, replacing  $\overline{z}$ with 
$  \zeta=\left(\zeta_{1},\zeta_{2},\dots,\zeta_{N}\right) $. Using $z'=\left(z_{1},z_{2},\dots,z_{N-1}\right)$ and $\zeta'=\left(\zeta_{1},\zeta_{2},\dots,\zeta_{N-1}\right),$ we consider the ring of the formal power series in $ x,y,t,r, z',z,\zeta$,  denoted by $\mathbb{C}\left[\left[   x,t,z',\zeta\right]\right]$, where $t,r$ is just the complexification of $x,y$. We assume $x,y,t,r$ near $0\in\mathbb{C}^{N}$ by an eventual translation.  We complexify the equations
(\ref{V21bu115}) and  (\ref{V22bu115}). We apply the Approximation Theorem of Artin\cite{A}. We obtain the desired convergence after few non-trivial computations.
 \section{Proof of Theorem \ref{t5}}
We  look  at  $z=\left(z_{1},z_{2},\dots,z_{N}\right)$ and $\overline{z}=\left(\overline{z}_{1},\overline{z}_{2},\dots,\overline{z}_{N}\right)$ as independent variables, replacing  $\overline{z}$ with 
$  \zeta=\left(\zeta_{1},\zeta_{2},\dots,\zeta_{N}\right) $. Using $z'=\left(z_{1},z_{2},\dots,z_{N-1}\right)$ and $\zeta'=\left(\zeta_{1},\zeta_{2},\dots,\zeta_{N-1}\right),$ we consider the ring of the formal power series in $ x,y,t,r, z',z,\zeta$,  denoted by $\mathbb{C}\left[\left[   x,t,z',\zeta\right]\right]$, where $t,r$ is just the complexification of $x,y$. We assume $x,y,t,r$ near $0\in\mathbb{C}^{N}$ by an eventual translation.  We complexify the equations
(\ref{V31}) and  (\ref{V32}). We apply the Approximation Theorem of Artin\cite{A}. We obtain the desired convergence after few non-trivial computations.
\section{Proof of Theorem \ref{t6}}
We  look  at  $z=\left(z_{1},z_{2},\dots,z_{N}\right)$ and $\overline{z}=\left(\overline{z}_{1},\overline{z}_{2},\dots,\overline{z}_{N}\right)$ as independent variables, replacing  $\overline{z}$ with 
$  \zeta=\left(\zeta_{1},\zeta_{2},\dots,\zeta_{N}\right) $. Using $z'=\left(z_{1},z_{2},\dots,z_{N-1}\right)$ and $\zeta'=\left(\zeta_{1},\zeta_{2},\dots,\zeta_{N-1}\right),$ we consider the ring of the formal power series in $ x,y,t,r, z',z,\zeta$,  denoted by $\mathbb{C}\left[\left[   x,t,z',\zeta\right]\right]$, where $t,r$ is just the complexification of $x,y$. We assume $x,y,t,r$ near $0\in\mathbb{C}^{N}$ by an eventual translation.  We complexify the equations
(\ref{Ve1}) and  (\ref{Ve2}). We apply the Approximation Theorem of Artin\cite{A}. We obtain the desired convergence after few non-trivial computations..
\section{Proof of Theorem \ref{t7}}
We  look  at  $z=\left(z_{1},z_{2},\dots,z_{N}\right)$ and $\overline{z}=\left(\overline{z}_{1},\overline{z}_{2},\dots,\overline{z}_{N}\right)$ as independent variables, replacing  $\overline{z}$ with 
$  \zeta=\left(\zeta_{1},\zeta_{2},\dots,\zeta_{N}\right) $. Using $z'=\left(z_{1},z_{2},\dots,z_{N-1}\right)$ and $\zeta'=\left(\zeta_{1},\zeta_{2},\dots,\zeta_{N-1}\right),$ we consider the ring of the formal power series in $ x,y,t,r, z',z,\zeta$,  denoted by $\mathbb{C}\left[\left[   x,t,z',\zeta\right]\right]$, where $t,r$ is just the complexification of $x,y$. We assume $x,y,t,r$ near $0\in\mathbb{C}^{N}$ by an eventual translation.  We complexify the equation
(\ref{V311a}). We apply the Approximation Theorem of Artin\cite{A}. We obtain the desired convergence after few non-trivial computations.
\section{Proof of Theorem \ref{t8}}
We  look  at  $z=\left(z_{1},z_{2},\dots,z_{N}\right)$ and $\overline{z}=\left(\overline{z}_{1},\overline{z}_{2},\dots,\overline{z}_{N}\right)$ as independent variables, replacing  $\overline{z}$ with 
$  \zeta=\left(\zeta_{1},\zeta_{2},\dots,\zeta_{N}\right) $. Using $z'=\left(z_{1},z_{2},\dots,z_{N-1}\right)$ and $\zeta'=\left(\zeta_{1},\zeta_{2},\dots,\zeta_{N-1}\right),$ we consider the ring of the formal power series in $ x,y,t,r, z',z,\zeta$,  denoted by $\mathbb{C}\left[\left[   x,t,z',\zeta\right]\right]$, where $t,r$ is just the complexification of $x,y$. We assume $x,y,t,r$ near $0\in\mathbb{C}^{N}$ by an eventual translation.  We complexify the equation
(\ref{Ve1ab}). We apply the Approximation Theorem of Artin\cite{A}. We obtain the desired convergence after few non-trivial computations.
\section{Proof of Theorem \ref{t9}}
We  look  at  $z=\left(z_{1},z_{2},\dots,z_{N}\right)$ and $\overline{z}=\left(\overline{z}_{1},\overline{z}_{2},\dots,\overline{z}_{N}\right)$ as independent variables, replacing  $\overline{z}$ with 
$  \zeta=\left(\zeta_{1},\zeta_{2},\dots,\zeta_{N}\right) $. Using $z'=\left(z_{1},z_{2},\dots,z_{N-1}\right)$ and $\zeta'=\left(\zeta_{1},\zeta_{2},\dots,\zeta_{N-1}\right),$ we consider the ring of the formal power series in $ x,y,t,r, z',z,\zeta$,  denoted by $\mathbb{C}\left[\left[   x,t,z',\zeta\right]\right]$, where $t,r$ is just the complexification of $x,y$. We assume $x,y,t,r$ near $0\in\mathbb{C}^{N}$ by an eventual translation.  We complexify the equation
(\ref{V311}) and  (\ref{V321}. We apply the Approximation Theorem of Artin\cite{A}. We obtain the desired convergence after few non-trivial computations.
\section{Proof of Theorem \ref{t10}}
We  look  at  $z=\left(z_{1},z_{2},\dots,z_{N}\right)$ and $\overline{z}=\left(\overline{z}_{1},\overline{z}_{2},\dots,\overline{z}_{N}\right)$ as independent variables, replacing  $\overline{z}$ with 
$  \zeta=\left(\zeta_{1},\zeta_{2},\dots,\zeta_{N}\right) $. Using $z'=\left(z_{1},z_{2},\dots,z_{N-1}\right)$ and $\zeta'=\left(\zeta_{1},\zeta_{2},\dots,\zeta_{N-1}\right),$ we consider the ring of the formal power series in $ x,y,t,r, z',z,\zeta$,  denoted by $\mathbb{C}\left[\left[   x,t,z',\zeta\right]\right]$, where $t,r$ is just the complexification of $x,y$. We assume $x,y,t,r$ near $0\in\mathbb{C}^{N}$ by an eventual translation.  We complexify the equation
(\ref{Ve1a}) and  (\ref{Ve2a}). We apply the Approximation Theorem of Artin\cite{A}. We obtain the desired convergence after few non-trivial computations.


\begin{thebibliography}{BER96b}
 \bibitem{A} {\bf Artin,~M.} --- On the solutions of analytic equations. {\em Inv. Math.} {\bf 5} (1968), nr. 4, $277-291$. 
   \bibitem{BEL1}{\bf Baouendi,~M.S.;   Linda Preiss Rothschild}---Geometric properties of mappings between hypersurfaces in complex space. {\em Journal of Diff. Geom.} {\bf 31}, (1990),   $473-499$.
    \bibitem{BEL2}{\bf Baouendi,~M.S.;   Linda Preiss Rothschild}---Mappings of real algebraic hypersurfaces. {\em Journal of the A.M.S.} {\bf 8}, (1995),   $997-1015$.
 \bibitem{BEL3}{\bf Baouendi,~M.S.; Ebenfelt,~.P; Linda Preiss Rothschild}---Convergence and Finite Jet Determination of Formal CR Mappings. {\em Journal of The A.M.S.} {\bf 13}, (2000), no. 4, $697-723$.
\bibitem{BREbook} {\bf Baouendi,~M.S.; Ebenfelt,~P.; Rothschild,~L.P.} ---{\em Real Submanifolds in Complex Space and Their Mappings.} Princeton
Math. Ser. {\bf 47}, Princeton Univ. Press, 1999.
\bibitem{BMR1}{\bf Baouendi,~M.S.;  Mir,~N.; Rothschild,~L.P.}--- Reflection ideals and mappings between generic submanifolds in complex space. {\em     J. Geom. Anal.}, 12 (4) (2002), pp. $543-580$.
\bibitem{Bha1} {\bf Bharali, ~G.} --- Surfaces with CR singularities that are locally polynomially convex, {\em Michigan Math. J.} {\bf 53} (2005), $429-445$.
\bibitem{Bha2} {\bf Bharali, ~G.} --- Polynomial approximation, local polynomial convexity, and degenerate CR singularities, {\em J. Funct. Anal.} {\bf 256} (2006), $351-368$.
\bibitem{Bha3} {\bf Bharali, ~G. }---Polynomial approximation, local polynomial convexity, and degenerate CR singularities II, {\em Internat. J. Math.} {\bf 22} (2011), $1721-1733$.
\bibitem{Bi} {\bf Bishop,~E.} --- Differentiable Manifolds In Complex
Euclidian Space. {\em Duke Math. J.} {\bf 32} (1965), no. 1,
$1-21$. 
\bibitem{V1} {\bf Burcea,~V.} --- A normal form for a real $2$-codimensional submanifold in $\mathbb{C}^{N+1}$ near a CR singularity. {\em  Adv. in Math.}
 {\bf 243} (2013), $262-295$.
\bibitem{V11} {\bf Burcea,~V.} --- On a family of analytic discs attached to a real submanifold $M\subset\mathbb{C}^{N+1}$, {\em Methods and Applications of Analysis} {\bf 20}, 1, (2013), $69-78$. 
\par \quad\quad\quad\quad\quad\quad\quad\quad\quad\quad\quad\quad\quad\quad\quad\quad\quad\quad\quad\quad  (with an Erratum submitted for publication)
 
  \bibitem{bu2} {\bf Burcea,~V.} --- Normal Forms and Degenerate CR Singularities. {\em  Complex Variables and Elliptic Equations} {\bf 61} (2016), 9, $1314-1333$.
\bibitem{Dol1} {\bf Dolbeault, ~P.} --- On Levi-flat hypersurfaces with given boundary in $\mathbb{C}^{n}$.  {\em Sci. China Ser. A} {\bf 51}, (2008), no. 4, $541-552$.
\bibitem{Dol3} {\bf Dolbeault,~P.} --- Boundaries of Levi-flat hypersurfaces: special hyperbolic points. {\em Ann. Polon. Math. } {\bf 106}, (2012), nr. 1, $145-170$. 
 \bibitem{DTZ2} {\bf Dolbeault, ~P.; Tomassini, ~G.; Zaitsev, ~D.} ---On Levi-flat hypersurfaces with prescribed boundary.    {\em Pure Appl. Math. Q.}, {\bf 6}, (2010), no. 3,   $725-753$. 
 \bibitem{DTZ1} {\bf Dolbeault, ~P.; Tomassini, ~G.; Zaitsev, ~D.} ---Boundary problem for Levi flat graphs. {\em Indiana Univ. Math. J.}, {\bf 60} (2011), no. 1, $161-170$.
\bibitem{G1} {\bf Gong,~X.} --- On the convergence of normalizations of real analytic surfaces near hyperbolic complex tangents. {\em Comment. Math. Helv.},  {\bf 69} (1994), no. 1, $549-574$.
\bibitem{G2} {\bf Gong,~X.} --- Existence of real analytic surfaces with hyperbolic complex tangent that are formally   but not holomorphically equivalent to quadrics. {\em Indiana Univ.
Math. J.} {\bf 53} (2004), no. 1, $83-95$.
 \bibitem{GL}{\bf  Gong,~X.; Lebl,~J.}--- Normal forms for CR singular codimension-two Levi-flat submanifolds.  {\em Pacific J. Math} 275 (2015), no. 1, $115-165$. 
\bibitem{HK}{\bf  Huang,~X.; Krantz,~S. }--- On a problem of Moser. {\em Duke Math. J.}, {\bf 78},  (1995), no. 1, $213-228$.  
 
 
\bibitem{HY1} {\bf Huang,~X.; Yin,~W.} --- A codimension two CR singular submanifold that is formally equivalent to a symmetric quadric. {\em Int. Math. Res. Notices}  (2009), no. 15, $2789-2828$. 
 \bibitem{HY2221} {\bf Huang,~X.; Yin,~W.} --- Flattening of CR singular points and the analyticity of the local hull of holomorphy I. {\em  Math. Ann.} {\bf 365} (2016), no. 1-2, $381-399$. 
  \bibitem{HY2222} {\bf Huang,~X.; Yin,~W.} --- Flattening of CR singular points and the analyticity of the local hull of holomorphy II. {\em Adv. Math.} {\bf 308} (2017), $1009-1073$.
  
   
\bibitem{RI} {\bf  Kossovskiy,~I. ; Shafikov, ~R.} ---Divergent CR-Equivalences and Meromorphic Differential Equations.{\em    Journal of European Math. Society},{\bf 18}, (2016), no.12, $2785-2819$.
 \bibitem{ko2} {\bf Kol\`{a}\v{r}, ~M.} --- Local symmetries  of finite type hypersurfaces  in $\mathbb{C}^{2}$, {\em Science in China A} {\bf 49} (2006),   $1633-1641$. 
\bibitem{ko3} {\bf Kol\`{a}\v{r}, ~M.} --- Local equivalence  of symmetric hypersurfaces  in $\mathbb{C}^{2}$, {\em Trans. of The A.M.S.} {\bf 362} (2010), no. 6,  $2833-2843$. 
   \bibitem{BEL5}{\bf Kolar,~M.}---The Catlin Multitype and Biholomorphic Equivalence of Models. {\em  International Mathematical Research Notices} {\bf 18}, (2010),   $3530 -3548$. 

\bibitem{Mer}{\bf  Merker,~J.}--- Convergence of formal invertible CR mappings between minimal holomorphically nondegenerate real analytic hypersurfaces. {\em Int. J. Math. Math. Sci.} {\bf 26}, (2001), no. 5, $281-302$.    

\bibitem{Mey} {\bf Meylan, F.} --- A reflection principle in complex space for a class of hypersurfaces and mappings, {\em Paciﬁc J. Math.}, {\bf 169 } (1995), $135-160$. 
\bibitem{MNZ}{\bf  Meylan,~F.; Mir,~N.; Zaitsev,~D.}--- Approximation and convergence of formal CR-mappings. {\em Int. Math. Res. Not.} (2003), no. 4, $211-242$.  
 \bibitem{Mir1} {\bf Mir,~N.}--- Formal biholomorphic maps of real analytic hypersurfaces.{\em Math. Res. Lett.} {\bf 7} (2000), no. 2-3, $343-359$.
 \bibitem{Mir2} {\bf Mir,~N.}--- On the convergence of formal mappings. {\em Comm. Anal. Geom.} {\bf 10} (2002), no. 1, $23-59$.
 \bibitem{MW} {\bf Moser,~ J.; Webster,~S.} --- Normal forms for real surfaces in
$\mathbb{C}^{2}$ near complex tangents and hyperbolic surface
transformations. {\em Acta Math.} {\bf 150} (1983), $255-296$.
\bibitem{Mos} {\bf  Moser,~J.} ---Analytic Surfaces in $\mathbb{C}^{2}$ and their local hull
of holomorphy. {\em Ann. Acad. Sci. Fenn. Ser. A.I. Math.} {\bf
10} (1985), 397-410. 
\bibitem{Sla2} {\bf Slapar,~M.} ---On Stein Neighborhood Basis of Real Surfaces. {\em Math. Z.}
 {\bf 247} (2004), no. $4$, $863-879$.
\bibitem{Su}{\bf  Suny\'{e},~J.C.}--- On formal maps between generic submanifolds in complex space. {\em J. Geom. Anal.} {\bf 19}, (2009), no. 4, $944-962$. 

 \bibitem{D11} {\bf Zaitsev,~D.} ---  Normal forms of non-integrable almost CR structures, {\em Amer. J. Math.} {\bf
134} (2012), no.4, $915-947$. 
\bibitem{D1} {\bf Zaitsev,~D.} --- New Normal Forms for Levi-nondegenerate
Hypersurfaces. {\em Several Complex Variables and Connections with
PDE Theory and Geometry}. Complex analysis-Trends in Math.,
  Birkhäuser/Springer Basel AG, Basel, pp. $321-340$,  (2010).
 \end{thebibliography}
\end{document}